\documentclass[12pt]{amsart}
\usepackage{amssymb}

\usepackage[mathcal]{euscript}

\newtheorem{thm}{Theorem}[section]
\newtheorem{prop}[thm]{Proposition}
\newtheorem{cor}[thm]{Corollary}

\newtheorem{mainthm}{Main Theorem}

\theoremstyle{definition}
\newtheorem{defn}[thm]{Definition}
\newtheorem{rem}[thm]{Remark}
\newtheorem{exmp}[thm]{Example}

\newtheorem*{notn}{Notation}

\numberwithin{equation}{section}

\renewcommand{\c}{\gamma}
\renewcommand{\d}{\delta}
\newcommand{\e}{\epsilon}
\newcommand{\f}{\varphi}
\renewcommand{\i}{^{-1}}

\newcommand{\C}{\mathcal{C}}
\newcommand{\D}{\mathcal{D}}
\newcommand{\E}{\mathcal{E}}
\newcommand{\F}{\mathbf{F}}

\newcommand{\Z}{\mathbf{Z}}

\renewcommand{\)}{\right)}
\renewcommand{\=}{\equiv}

\newcommand{\adt}{\operatorname{adt}}
\newcommand{\Aut}{\operatorname{Aut}}
\newcommand{\rad}{\operatorname{rad}}

\newcommand{\dividing}{\mid}
\newcommand{\gen}[1]{\left\langle #1 \right\rangle}

\newcommand{\intersect}{\cap}
\newcommand{\isom}{\cong}

\newcommand{\noqed}{\renewcommand{\qed}{}}
\newcommand{\normal}{\vartriangleleft}

\newcommand{\set}[1]{\left\{ {#1} \right\}}
\newcommand{\tr}{\operatorname{trace}}
\newcommand{\wt}[1]{\left| {#1} \right|}

\newcommand{\sdcp}{%
\,\mbox{{\footnotesize $\bigcirc$}\hspace*{-.715em}\raisebox{.05em}{\tiny
$\pm$}}\mspace{4.9mu}}
\newcommand{\cp}{\,\mbox{{\footnotesize $\bigcirc$}}\,}

\begin{document}

%%%%% intro, background

\title[Class 2 and small Frattini Moufang loops]{Class 2 Moufang
	loops, small Frattini Moufang loops, and code loops}
\author{Tim Hsu}
\address{\hskip-\parindent Department of Mathematics \\ University of Michigan \\
	Ann Arbor, MI 48109}
\email{timhsu@@math.lsa.umich.edu}

\subjclass{Primary 20N05; Secondary 20D15, 20D08}
\keywords{Centrally nilpotent Moufang loops of class 2, Moufang
	$p$-loops, small Frattini Moufang loops, code loops}

\date{\today}

\begin{abstract}
Let $L$ be a Moufang loop which is centrally nilpotent of class 2.  We
first show that the nuclearly-derived subloop (normal associator
subloop) $L^*$ of $L$ has exponent dividing 6.  It follows that $L_p$
(the subloop of $L$ of elements of $p$-power order) is associative for
$p>3$.  Next, a loop $L$ is said to be a {\it small Frattini Moufang
loop}, or SFML, if $L$ has a central subgroup $Z$ of order $p$ such
that $C\isom L/Z$ is an elementary abelian $p$-group.  $C$ is thus
given the structure of what we call a {\it coded vector space}, or
CVS.  (In the associative/group case, CVS's are either orthogonal
spaces, for $p=2$, or symplectic spaces with attached linear forms,
for $p>2$.)  Our principal result is that every CVS may be obtained
from an SFML in this way, and two SFML's are isomorphic in a manner
preserving the central subgroup $Z$ if and only if their CVS's are
isomorphic up to scalar multiple.  Consequently, we obtain the fact
that every SFM 2-loop is a code loop, in the sense of Griess, and we
also obtain a relatively explicit characterization of isotopy in SFM
3-loops.  (This characterization of isotopy is easily extended to
Moufang loops of class 2 and exponent 3.)  Finally, we sketch a method
for constructing any finite Moufang loop which is centrally nilpotent
of class 2.
\end{abstract}

\maketitle

\section{Introduction}
\label{sect:intro}

The loops (groups without associativity) characterized by the
near-associativity property
\begin{equation}
(xy)(zx) = x((yz)x)
\end{equation}
are known as {\it Moufang loops} (see Pflugfelder~\cite[Ch.\
IV]{hop:loops}).  Many aspects of group theory may be generalized to
Moufang loops, and among these aspects is the theory of {\it central
nilpotence\/} (Bruck~\cite[Ch.\ VI]{rhb:loops}), the loop
generalization of nilpotence in groups.  Centrally nilpotent Moufang
loops have been studied often, and provide many of the basic examples
of finite Moufang loops.  (See, for instance, Bruck~\cite[Ch.\ VIII,
Thm.\ 10.1]{rhb:loops}, Chein~\cite[II.4]{oc:exloops},
Pflugfelder~\cite[Ch.\ IV]{hop:loops}, and Smith~\cite[p.\
181]{jdhs:quasichar}.)  Of particular relevance to this paper is the
work of Glauberman and Wright~\cite{gg:oddloopsII,ggcrbw:nilploops},
who extended many of the standard theorems about finite nilpotent
groups to finite centrally nilpotent Moufang loops.

One class of centrally nilpotent Moufang loops which has particularly
interesting applications to finite group theory is the class of {\it
code loops}.  The first code loop to be recognized as such was the
{\it Parker loop\/} (named after its discoverer R. A. Parker), which
played a key role in Conway's construction of the Monster finite
simple group~\cite{jhc:monster}.  Subsequently,
Griess~\cite{rlg:codeloops} defined code loops to be certain central
extensions of doubly even codes, providing the first published proof
of their existence, and then went on~\cite{rlg:codetwolocal} to use
code loops to construct 2-local subgroups of several other sporadic
groups.  For more on code loops and finite groups, including further
references, see Griess~\cite{rlg:codetwolocal} and
Richardson~\cite{tmr:oddcodeloops}.

In this paper, we present some new results on centrally nilpotent
Moufang loops of class 2, and apply these results by generalizing the
theory of extraspecial groups (see, for instance,
Aschbacher~\cite[Ch.\ 8]{ma:finite}) to {\it small Frattini Moufang
loops}.  For the convenience of the reader, we now summarize our main
results.  (The reader who is unfamiliar with the notation and
terminology used here may first wish to read Section
\ref{sect:background}.)

Let $L$ be a Moufang loop which is centrally nilpotent of class 2,
that is, a Moufang loop $L$ such that the quotient of $L$ by its
center $Z(L)$ is an abelian group; and let $L_p$ be the set of all
elements of $L$ whose order is a power of $p$.  Recall that the
nuclearly-derived subloop, or normal associator subloop, of $L$, which
we denote by $L^*$, is the smallest normal subloop of $L$ such that
$L/L^*$ is associative (is a group).  Recall also that the torsion
subloop (subloop of finite order elements) of $L$ is isomorphic to the
(restricted) direct product of the subloops $L_p$, where $p$ runs over
all primes (Thm.\ 6.2 of Bruck~\cite{rhb:somemoufang}, our Theorem
\ref{thm:torsionsubloop}, or in the finite case, Cor.\ 1 of Glauberman
and Wright~\cite{ggcrbw:nilploops}).

In Section \ref{sect:classtwo}, we show that the commutator (resp.\
associator) function is a ``symplectic'' and ``multilinear'' function
on $L/Z(L)\times L/Z(L)$ (resp.\ $L/Z(L)\times L/Z(L)\times L/Z(L)$)
(Theorem \ref{thm:sympmulti}).  As a consequence, we have:

\begin{mainthm}\label{mainthm:assocsix}
Let $L$ be a Moufang loop which is centrally nilpotent of class 2.
Then $L^*$ (as defined above) has exponent dividing 6.  In particular,
for $p>3$, $L_p$ (as defined above) is associative (is a group).
\end{mainthm}

Compare the result of Bruck~\cite[VIII.2]{rhb:loops} that the cube of
every associator of a commutative Moufang loop is trivial.  (In fact,
to prove Main Theorem \ref{mainthm:assocsix}, we use another case of
the same formulas Bruck used to obtain that result.)  We also note
that Main Theorem \ref{mainthm:assocsix} is, in some sense, the best
possible result of this type, since Example 3 of VII.5 of
Bruck~\cite{rhb:loops} gives a construction of nonassociative finite
Moufang $p$-loops of class 3 for all $p>2$.

In Sections \ref{sect:eml}--\ref{sect:codedextensions}, we focus on
{\it small Frattini Moufang loops\/} (also known as SFM loops, or
SFML's), which are Moufang $p$-loops $L$ with a central subgroup $Z$
of order $p$ such that $C\isom L/Z$ is an elementary abelian group.
(Note that we often think of $Z$ and $C$ as part of the structure of
$L$.)  SFML's are a class of Moufang loops often found ``in nature.''
For instance, every extraspecial Moufang loop (Definition
\ref{defn:eml}) is an SFML.  Also, every code loop is an SFM 2-loop,
and conversely:

\begin{mainthm}\label{mainthm:evenSFMLiscodeloop}
Every SFM 2-loop is isomorphic to a code loop.
\end{mainthm}

Compare Thm.\ 14 of Griess~\cite{rlg:codeloops}, which shows that
every loop constructed by ``Parker's procedure'' (Defn.\ 13 of
Griess~\cite{rlg:codeloops}) is isomorphic to a code loop.  For a more
detailed comparison, see Remarks~\ref{rem:griesspf} and
\ref{rem:codesize}.

The key to Main Theorem \ref{mainthm:evenSFMLiscodeloop}, and also to
Main Theorem \ref{mainthm:CVSisSFML}, below, is the fact that $L$
gives $C$ the structure of a {\it coded vector space}, or {\it CVS},
over $\F_p$.  A CVS is a 4-tuple $(C,\sigma,\chi,\alpha)$, where $C$
is a vector space over $\F_p$, and $\sigma$, $\chi$, and $\alpha$ are
``symplectic 1-, 2-, and 3-forms'' on $C$ which are either
multilinear, for $p>2$, or related by polarization, for $p=2$.

The relationship between SFML's and CVS's can be stated as follows.

\begin{mainthm}\label{mainthm:CVSisSFML}
Every CVS can be obtained from some SFML in the manner described
above.  Furthermore, let $L$ and $M$ be SFML's, with distinguished
central subgroups $Z_L$ and $Z_M$, and associated CVS's $C_L$ and
$C_M$.  Then there is an isomorphism $\phi: L \rightarrow M$ such that
$\phi(Z_L)=Z_M$ if and only if $C_L$ and $C_M$ are isomorphic up to
scalar multiple (action of $\Aut(Z_L)=\Aut(Z_M)$).
\end{mainthm}

The proofs of Main Theorems \ref{mainthm:evenSFMLiscodeloop} and
\ref{mainthm:CVSisSFML} may be summarized as follows.
\begin{list}{}{%
	\setlength{\leftmargin}{8em}
	\setlength{\labelwidth}{5em}
	\setlength{\labelsep}{\leftmargin}
	\addtolength{\labelsep}{-\labelwidth}
	\addtolength{\labelsep}{-\parindent}
}
\item[\bf Section \ref{sect:eml}.]  We define the notion of a {\it
coded extension\/} of a CVS (Definition \ref{defn:codedext}), and show
that every SFML is a coded extension of a CVS, and vice versa (Theorem
\ref{thm:codedextisSFML}).
\item[\bf Section \ref{sect:cvsiscode}.]  We show that every CVS over
$\F_2$ can be obtained from a doubly even code (Theorems
\ref{thm:codeisCVS} and \ref{thm:CVSiscode}).  Main Theorem
\ref{mainthm:evenSFMLiscodeloop} follows.
\item[\bf Section \ref{sect:codedextensions}.]  We show that every CVS
has a unique coded extension (Theorems \ref{thm:codedextunique} and
\ref{thm:codedextexist}).  Main Theorem \ref{mainthm:CVSisSFML}
follows.
\end{list}

We remark that the main technical tool used in Section
\ref{sect:codedextensions}, the {\it semidirect central product}, is
also useful for doing calculations in SFML's, especially code loops.
In~\cite{th:codeloops}, we will address the general topic of
decompositions of SFML's as semidirect central products of groups.  In
particular, we will give some short explicit constructions of code
loops, including a Turyn-type construction for the Parker loop.  For
more details, see Remark~\ref{rem:sdcp}.

Now, as mentioned above, every SFM 3-loop $L$ is a coded extension of
some CVS over $\F_3$, say, $(C,\sigma(c),\chi(c,d),\alpha(c,d,e))$.
For any $k\in C$, we define the {\it adjoint translate\/} $\adt_k(C)$
of $C$ to be the CVS
$(C,\sigma(c),\chi(c,d)+\alpha(c,k,d),\alpha(c,d,e))$.  A
straightforward application of Main Theorem~\ref{mainthm:CVSisSFML}
then gives the following characterization of isotopy in SFM 3-loops
(Section \ref{sect:isotopy}).

\begin{mainthm}\label{mainthm:isotopy}
Let $L$ be a coded extension of a CVS $(C,\sigma,\chi,\alpha)$ over
$\F_3$.  Then up to isomorphism, the loop-isotopes of $L$ are
precisely the coded extensions of the adjoint translates of $C$.  In
particular, $\sigma$ and $\alpha$ are ``isotopy invariants'' of $L$.
\end{mainthm}

Furthermore, since Main Theorem~\ref{mainthm:CVSisSFML} may be
generalized directly to Moufang loops of class 2 and exponent 3, Main
Theorem \ref{mainthm:isotopy} may also be extended to this situation.
See Theorem \ref{thm:expthreeisotopy} for a precise statement.

We conclude in Section \ref{sect:constructclasstwo} by giving a
construction which can be used to obtain any finite Moufang loop of
class 2.  This construction generalizes much previous work, mostly in
the commutative case.  For instance, compare
B\'{e}n\'{e}teau~\cite[IV.3]{lb:cml}, Bruck~\cite{rhb:contloops},
Chein~\cite[II.4]{oc:exloops}, and Ray-Chaudhuri and
Roth~\cite{dkrcrr:hallcml}.

\section{Background and notation}
\label{sect:background}

First, we set some conventions and notation to be used throughout.

\begin{notn}
Let $p$ be a prime.  $\F_p$ denotes the field of order $p$, and
$\F_p^\times$ its nonzero elements.  Following group-theoretic custom,
unless otherwise specified, we think of $\F_p$ as the group of order
$p$, and $\F_p^\times$ as the automorphism group of $\F_p$.  In this
context, we identify the vector space $\F_p^k$ with the elementary
abelian $p$-group of rank $k$, and we write vector addition in
$\F_p^k$ multiplicatively, with the zero vector written as 1.

If $a,b,c,\dots$ are elements or subsets of a loop (resp.\ vector
space), $\gen{a,b,c,\dots}$ denotes the subloop (resp.\ subspace)
generated by $a,b,c,\dots$.
\end{notn}

For those less familiar with loop theory, and for the purpose of
establishing notation and terminology, we also review some definitions
and results in loop theory, using Pflugfelder~\cite{hop:loops},
Bruck~\cite{rhb:loops}, and Chein, Pflugfelder, and
Smith~\cite{cps:loopsappl} as our standard sources.

\begin{defn}
An {\it inverse property loop}, or in this paper, simply a {\it loop},
is a set $L$ with a binary operation (written as juxtaposition) having
a unique identity element and unique two-sided inverses.  (Note that
the term ``inverse'' means that $a\i(ax)=(xa)a\i=x$.)  In other words,
a loop is a group minus associativity.
\end{defn}

Many concepts of group theory may be generalized to loop theory; we
highlight the following ones.

\begin{defn}
For loop elements $\c$, $\d$, and $\e$, we define the {\it
commutator\/} $[\c,\d]$ to be $(\d\c)\i(\c\d)$ and the {\it
associator\/} $[\c,\d,\e]$ to be $(\c(\d\e))\i((\c\d)\e)$.  In other
words,
\begin{align}
\c\d &= (\d\c)[\c,\d], \\
(\c\d)\e &= (\c(\d\e))[\c,\d,\e], \\
\c(\d\e) &= ((\c\d)\e)[\c,\d,\e]\i.
\end{align}
(The inexperienced reader should note the inverse in the last
formula.)
\end{defn}

\begin{defn}\label{defn:centnuc}
Let $L$ be a loop.  The {\it nucleus\/} of $L$ (denoted by $N(L)$) is
the set of all $z\in L$ such that $[z,x,y]=[x,z,y]=[x,y,z]=1$ for all
$x,y\in L$; and the {\it center\/} of $L$ (denoted by $Z(L)$) is
defined to be the set of all $z\in N(L)$ such that $[z,x]=1$ for all
$x,y\in L$.
\end{defn}

If $L$ is a loop, it can be shown (see
Pflugfelder~\cite[I.3]{hop:loops}) that $N(L)$ is a subgroup of $L$,
and that $Z(L)$ is an abelian subgroup of $N(L)$.

\begin{defn}
A {\it normal\/} subloop of a loop $L$ is any subloop of $L$ which is
the kernel of some homomorphism from $L$ to a loop.
\end{defn}

For instance, any {\it central subgroup\/} (subgroup of $Z(L)$) of a
loop $L$ is normal in $L$ (Pflugfleder~\cite[I.7]{hop:loops}).

\begin{defn}\label{defn:derived}
Let $L$ be a loop.  We define the {\it centrally-derived subloop\/}
(or normal commutator-associator subloop) of $L$ to be the smallest
normal subloop $L'\normal L$ such that $L/L'$ is an abelian group.
Similarly, we define the {\it nuclearly-derived subloop\/} (or normal
associator subloop) of $L$ to be the smallest normal subloop
$L^*\normal L$ such that $L/L^*$ is associative (is a group).
\end{defn}

See Bruck~\cite[Ch.\ VI]{rhb:loops} for a proof that $L'$ and $L^*$
are well-defined.  Note that it follows from the isomorphism theorems
for loops (see Pflugfelder~\cite[I.7]{hop:loops}) that $L'$ (resp.\
$L^*$) is the smallest normal subloop of $L$ containing all $[\c,\d]$
and $[\c,\d,\e]$ (resp.\ all $[\c,\d,\e]$), where $\c,\d,\e$ run over
all elements of $L$.

We will use Bruck's theory of central nilpotence~\cite[Ch.\
VI]{rhb:loops}, as described in the following definitions and Theorem
\ref{thm:frattini}.

\begin{defn}
Let $L$ be a loop.  The {\it upper central series\/} $\set{Z_i}$ of
$L$ is defined by letting $Z_{0}=1$ and letting $Z_{i+i}$ be the
unique subloop of $L$ containing $Z_{i}$ such that
$Z_{i+i}/Z_{i}=Z(L/Z_{i})$.  We say that $L$ is {\it centrally
nilpotent of class $n$}, or simply of class $n$, if there exists $n$
such that $Z_{n}=L$ and $Z_{n-1}\neq L$.
\end{defn}

For instance, $L$ is of class 2 if and only if $L/Z(L)$ is an abelian
group and $L$ is not, that is, if and only if $1<L'\le Z(L)$.

\begin{defn}
The {\it Frattini subloop\/} $\Phi(L)$ of a loop $L$ is defined to be
the set of {\it non-generators\/} of $L$, that is, the set of all
$x\in L$ such that for any subset $S$ of $L$, $L=\gen{x,S}$ implies
$L=\gen{S}$.
\end{defn}

\begin{thm}\label{thm:frattini}
Let $L$ be a finite centrally nilpotent loop.  Then $\Phi(L)\normal
L$, and $L/\Phi(L)$ is isomorphic to a subgroup of the direct product
of groups of prime order.
\end{thm}
\begin{proof}
This follows immediately from Thms.\ 2.1 and 2.2 of Ch.\ VI of
Bruck~\cite{rhb:loops}.
\end{proof}

We are particularly interested in loops of the following type.

\begin{defn}\label{defn:Moufang}
A loop $L$ is said to be {\it Moufang\/} if any, and therefore all (see
Pflugfelder~\cite[Ch.\ IV]{hop:loops}), of the following identities
hold for all $\c,\d,\e\in L$:
\begin{gather}
\begin{align}
((\d\c)\e)\c &= \d(\c(\e\c)), \label{eq:firstMoufang} \\
((\c\d)\c)\e &= \c(\d(\c\e)), \label{eq:Mkone}
\end{align} \\
(\c(\d\e))\c = (\c\d)(\e\c) = \c((\d\e)\c). \label{eq:lastMoufang}
\end{gather}
\end{defn}

\begin{defn}\label{defn:Moufangcenter}
The {\it Moufang center\/} of a Moufang loop $L$, denoted by $C(L)$,
is defined to be the set of all $z\in L$ such that $[z,x]=1$ for all
$x\in L$.
\end{defn}

Let $L$ be a Moufang loop.  Clearly, $Z(L)=N(L)\cap C(L)$.
Furthermore, it can also be shown (see Pflugfelder~\cite[Thm.\
IV.3.10]{hop:loops}) that $C(L)$ is a subloop of $L$.

Moufang loops have many near-associativity properties, such as the
following consequence of {\it Moufang's theorem\/} (see
Pflugfelder~\cite[Ch.\ IV]{hop:loops}).

\begin{thm}
Let $L$ be a Moufang loop.  Then $L$ is di-associative; that is, for
$x,y\in L$, $\gen{x,y}$ is associative.  In particular, $L$ is
power-associative; that is, $x^n$ is well-defined.  \qed
\end{thm}

We will also use the Lagrangian property of di-associative loops
(Bruck~\cite[Thm.\ V.1.2]{rhb:loops}), stated as:

\begin{thm}\label{thm:Lagrange}
Let $L$ be a finite di-associative (e.g., Moufang) loop.  Then the
order of any element of $L$ divides the order of $L$.  \qed
\end{thm}

\begin{defn}
Let $L$ be a power-associative loop, and let $p$ be a prime.  We say
that $L$ is a {\it $p$-loop\/} if every element of $L$ has order a
power of $p$.
\end{defn}

Let $L$ be a di-associative loop.  Because of Theorem
\ref{thm:Lagrange}, if $L$ has order a power of $p$, then $L$ is a
finite $p$-loop.  Conversely, if $L$ is a finite centrally nilpotent
$p$-loop, the isomorphism theorems for loops imply that the order of
$L$ is a power of $p$.  Therefore, since most of the loops we consider
are finite centrally nilpotent Moufang loops, we will usually treat
the concepts of having order a power of $p$ and being a finite
$p$-loop as interchangable.  (In fact, all finite Moufang $p$-loops
are centrally nilpotent; see Glauberman~\cite[Thm.\ 4]{gg:oddloopsII}
and Glauberman and Wright~\cite{ggcrbw:nilploops}.)

Finally, we define one last important concept of loop theory.

\begin{defn}
A triple $(U,V,W)$ of bijections from a loop $L$ to a loop $M$ (whose
operation is denoted by $\circ$) is called an {\it isotopism\/} if,
for all $x,y\in L$, $(xU)\circ(yV)=(xy)W$.  If an isotopism exists
from $L$ to $M$, we say that $M$ is an {\it isotope\/} of $L$, or that
$L$ and $M$ are {\it isotopic}.
\end{defn}

It is worth noting that isotopy plays no role in group theory because
every loop-isotope of a group $G$ is isomorphic to $G$ (see
Pflugfelder~\cite[Cor.\ III.2.3]{hop:loops}).  More generally, any
loop which is isomorphic to all of its loop-isotopes is called a {\it
$G$-loop}.  See Pflugfelder~\cite[Ch.\ III]{hop:loops} for more on
isotopy.

\section{Moufang loops of class 2}
\label{sect:classtwo}

We first quote the following result, due to Bruck.

\begin{prop}\label{prop:rhblem}
Let $L$ be a Moufang loop such that $[[\c,\d,\e],\c]=1$ for all
$\c,\d,\e\in L$.  Then for all $\c,\d,\e\in L$, $[\c,\d,\e]$ is
central in $\gen{\c,\d,\e}$, and the following identities hold for all
$n\in\Z$:
\begin{gather}
[\c,\d,\e] = [\d,\e,\c] = [\d,\c,\e]\i \label{eq:assocskew} \\
[\c^n,\d,\e] = [\c,\d,\e]^n \label{eq:assocpowerlin} \\
[\c\d,\e] = [\c,\e] [[\c,\e],\d] [\d,\e] [\c,\d,\e]^3
	\label{eq:commmultilin}
\end{gather}
\end{prop}
Note that part of the statement of \eqref{eq:commmultilin} is that the
right hand side gives the same result, no matter how the terms are
associated.
\begin{proof}
This follows from Lemma VII.5.5 of Bruck~\cite{rhb:loops}.
\end{proof}

For the rest of this section, let $L$ be a Moufang loop with a fixed
central subgroup $Z$ such that $C\isom L/Z$ is an abelian group.
Clearly, such a loop is centrally nilpotent of class 2, and
conversely, for any Moufang loop $L$ of class 2, we may take $Z=Z(L)$.
By convention, the letters $\c$, $\d$, $\e$, and $\f$ refer to
elements of $L$, and their images in the quotient $C$ are denoted by
$c$, $d$, $e$, and $f$, respectively.

\begin{defn}\label{defn:chialpha}
We define functions $\chi: C \times C \rightarrow Z$ and $\alpha:
C\times C \times C\rightarrow Z$ by the following formulas.
\begin{align}
\chi(c,d) &= [\c,\d], \\
\alpha(c,d,e) &=[\c,\d,\e].
\end{align}
Note that $\chi$ and $\alpha$ are well-defined because for any $z\in
Z$, $[\c z,\d,\e]=[\c,\d,\e]$, and so on.
\end{defn}

The following key theorem says that the functions $\chi$ and
$\alpha$ are ``symplectic'' (\eqref{eq:chisymp} and
\eqref{eq:alphasymp}), ``skew-symmetric'' (\eqref{eq:chiskew} and
\eqref{eq:alphaskew}), ``power-multilinear'' (\eqref{eq:chipowerlin} and
\eqref{eq:alphapowerlin}), and
related by ``polarization'' \eqref{eq:chimultilin}; and that $\alpha$
is multilinear \eqref{eq:alphamultilin}.

\begin{thm}\label{thm:sympmulti}
For all $c,d,e,f \in C$ and all $n\in\Z$, we have:
\begin{align}
\chi(c,c) &= 1, \label{eq:chisymp} \\
\chi(c,d) &= \chi(d,c)\i, \label{eq:chiskew} \\
\chi(c^n,d) &= \chi(c,d)^n, \label{eq:chipowerlin} \\
\chi(cd,e) &= \chi(c,e)\chi(d,e)\alpha(c,d,e)^3,
	\label{eq:chimultilin}
\end{align}
and
\begin{align}
\alpha(c,d,d) &= \alpha(d,c,d) = \alpha(d,d,c) = 1,
	\label{eq:alphasymp} \\
\alpha(c,d,e) &= \alpha(d,c,e)\i = \alpha(d,e,c),
	\label{eq:alphaskew} \\
\alpha(c^n,d,e) &= \alpha(c,d,e)^n, \label{eq:alphapowerlin} \\
\alpha(cd,e,f) &= \alpha(c,e,f)\alpha(d,e,f). \label{eq:alphamultilin}
\end{align}
\end{thm}
\begin{proof}
We first note that \eqref{eq:chisymp} and \eqref{eq:chiskew} are easy,
\eqref{eq:chimultilin} follows from \eqref{eq:commmultilin} and the
fact that $L'$ is central, \eqref{eq:alphasymp} follows from
di-associativity, \eqref{eq:alphaskew} follows from
\eqref{eq:assocskew}, and \eqref{eq:alphapowerlin} follows from
\eqref{eq:assocpowerlin}.  Furthermore, \eqref{eq:chimultilin},
\eqref{eq:alphasymp}, and \eqref{eq:alphapowerlin} imply
\begin{equation}\label{eq:chipowerlintrick}
\begin{split}
\chi(c^{n+1},d) &= \chi(c^n,d)\chi(c,d)\alpha(c^n,c,d)^3 \\
	&= \chi(c^n,d)\chi(c,d),
\end{split}
\end{equation}
so \eqref{eq:chipowerlin} follows by induction on positive and
negative $n$.

It remains to prove \eqref{eq:alphamultilin}.  Now, by definition,
\begin{equation}
\alpha(cd,e,f) = ((\c\d)(\e\f))\i (((\c\d)\e)\f),
\end{equation}
and
\begin{equation}\label{eq:assocerror}
\begin{split}
((\c\d)\e)\f &= (\c(\d\e))\f \cdot \alpha(c,d,e) \\
	&= \c((\d\e)\f) \cdot \alpha(c,d,e)\alpha(c,de,f) \\
	&= \c(\d(\e\f)) \cdot
		\alpha(c,d,e)\alpha(c,de,f)\alpha(d,e,f) \\
	&= (\c\d)(\e\f) \cdot
		\alpha(c,d,e)\alpha(c,de,f)\alpha(d,e,f)\alpha(c,d,ef)\i,
\end{split}
\end{equation}
which means that
\begin{equation}
\alpha(cd,e,f) = \alpha(c,d,e)\alpha(c,de,f)
	\alpha(d,e,f)\alpha(c,d,ef)\i. \label{eq:fivecycle}
\end{equation}
We claim that \eqref{eq:alphamultilin} is a consequence of
\eqref{eq:fivecycle}.

To prove this claim, by substituting first $c=w$, $d=x$, $e=y$, and
$f=z$, and then $c=x$, $d=y$, $e=z$, and $f=w$, into
\eqref{eq:fivecycle}, we get
\begin{align}
\alpha(wx,y,z) &= \alpha(w,x,y)\alpha(w,xy,z)
	\alpha(x,y,z)\alpha(w,x,yz)\i, \label{eq:fivesub} \\
\alpha(xy,z,w) &= \alpha(x,y,z)\alpha(x,yz,w)
	\alpha(y,z,w)\alpha(x,y,zw)\i. \label{eq:fivesubperm}
\end{align}
Since skew-symmetry implies $\alpha(w,xy,z)=\alpha(xy,z,w)$, we may
substitute the right-hand side of \eqref{eq:fivesubperm} for the
second term in the right-hand side of \eqref{eq:fivesub}.  Applying
skew-symmetry to collect terms, we obtain
\begin{equation}
\alpha(wx,y,z) = \alpha(wz,y,x)
	\alpha(w,x,y)\alpha(w,y,z)\alpha(x,y,z)^2.
	\label{eq:exchange}
\end{equation}
We call \eqref{eq:exchange} the {\it exchange identity}, since it
implies that we may exchange the $x$ and the $z$ in $\alpha(wx,y,z)$
at the cost of adding the other terms on the right-hand side of
\eqref{eq:exchange}.

Using the exchange identity and skew-symmetry, we see that
\begin{equation} \label{eq:killde}
\begin{split}
\alpha(c,de,f) &= \alpha(de,f,c) \\
	&= \alpha(dc,f,e)\alpha(d,e,f)\alpha(d,f,c)\alpha(e,f,c)^2 \\
	&= \alpha(cd,e,f)\i\alpha(d,e,f)\alpha(c,d,f)\alpha(c,e,f)^2.
\end{split}
\end{equation}
Applying exchange and skew-symmetry again, we have
\begin{equation}\label{eq:prekillef}
\begin{split}
\alpha(c,d,ef)\i &= \alpha(ef,d,c) \\
	&= \alpha(ec,d,f)\alpha(e,f,d)\alpha(e,d,c)\alpha(f,d,c)^2 \\
	&= \alpha(ce,f,d)\i
		\alpha(d,e,f)\alpha(c,d,e)\i\alpha(c,d,f)^{-2},
\end{split}
\end{equation}
and applying exchange and skew-symmetry to the first term of the last
expression in \eqref{eq:prekillef}, we have
\begin{equation}\label{eq:killef}
\begin{split}
\alpha(c,d,ef)\i &= \alpha(cd,f,e)\i \alpha(c,e,f)\i
		\alpha(c,f,d)\i \alpha(e,f,d)^{-2} \\
	&\qquad
		\cdot\alpha(d,e,f)\alpha(c,d,e)\i
		\alpha(c,d,f)^{-2}, \\
	&= \alpha(cd,e,f) \alpha(c,d,e)\i \alpha(c,d,f)\i
		\alpha(c,e,f)\i \alpha(d,e,f)\i.
\end{split}
\end{equation}

Finally, substituting \eqref{eq:killde} and \eqref{eq:killef} into
\eqref{eq:fivecycle}, we get
\begin{equation}
\begin{split}
\alpha(cd,e,f) &= \alpha(c,d,e) \\
&\qquad
	\cdot\alpha(cd,e,f)\i\alpha(d,e,f)\alpha(c,d,f)
	\alpha(c,e,f)^2 \\
&\qquad
	\cdot\alpha(d,e,f) \\
&\qquad
	\cdot\alpha(cd,e,f)
	\alpha(c,d,e)\i \alpha(c,d,f)\i
	\alpha(c,e,f)\i \alpha(d,e,f)\i \\
&= \alpha(c,e,f)\alpha(d,e,f),
\end{split}
\end{equation}
and the theorem follows.
\end{proof}

\begin{rem}
Note that \eqref{eq:alphamultilin}, which is really the only new
formula in Theorem~\ref{thm:sympmulti}, has been previously obtained
in several special cases, such as the commutative case
(Smith~\cite{jdhs:associators}) and the code loop case
(Griess~\cite[Lem.\ 15]{rlg:codeloops}).
\end{rem}

\begin{rem}
It may be instructive to consider the following method of proving
Theorem~\ref{thm:sympmulti} without relying on Bruck's formulas
(Proposition \ref{prop:rhblem}).  (In fact, this is how the author
first discovered Theorem~\ref{thm:sympmulti}.)  Now, given the
skew-symmetry of $\alpha$ (eq.\ \eqref{eq:alphaskew}), we can obtain
\eqref{eq:alphamultilin} as above, and we can obtain
\eqref{eq:chimultilin} by calculating the $\chi$ and $\alpha$ terms
needed to change $(\c\d)\e$ to $\e(\c\d)$, as we did in
\eqref{eq:assocerror}.  Furthermore, \eqref{eq:chipowerlin} and
\eqref{eq:alphapowerlin} follow from \eqref{eq:chimultilin} and
\eqref{eq:alphamultilin} and di-associativity.  Therefore, the crux of
the proof lies in obtaining \eqref{eq:alphaskew}.  In fact, it is here
that the Moufang property seems to be used most strongly, as the only
proofs of skew-symmetry of which the author is aware rely on the fact
that every inner mapping of a Moufang loop is a semi-endomorphism.
(This approach involves imitating one part of the proof of Moufang's
theorem; see, for example, Pflugfelder~\cite[IV.2.3]{hop:loops}.)
\end{rem}

\begin{rem}\label{rem:schneps}
Schneps (personal communication) has observed that
\eqref{eq:fivecycle} is yet another version of the ``pentagonal''
relation from monoidal categories, and that \eqref{eq:chimultilin} is
a version of the ``hexagonal'' relation from symmetric monoidal
categories.  See MacLane~\cite[Ch.\ VII]{sm:categories} for more on
these relations; see also Remark \ref{rem:conway}.
\end{rem}

In the rest of this section, we describe some of the consequences of
Theorem~\ref{thm:sympmulti}.

\begin{thm}\label{thm:gcd}
For $c,d,e\in C$ such that $c^k=d^m=e^n=1$, the order of $\chi(c,d)$
divides $\gcd(k,m)$, and the order of $\alpha(c,d,e)$ divides
$\gcd(k,m,n)$.
\end{thm}
\begin{proof}
From \eqref{eq:chipowerlin}, we have
\begin{equation}
\chi(c,d)^k = \chi(c^k,d) = \chi(1,d) = 1,
\end{equation}
and our commutator claim follows from skew-symmetry.  The same proof
works for our associator claim.
\end{proof}

For a prime $p$, define $L_p$ to be the set of all $x\in L$ such that
the order of $x$ is a power of $p$.  In the next two theorems (Theorems
\ref{thm:Lpsubloop} and \ref{thm:torsionsubloop}) we recover the class
2 case of results of Bruck~\cite[Thm.\ 6.2]{rhb:somemoufang} and
Glauberman and Wright~\cite[Cor.\ 1]{ggcrbw:nilploops}.

\begin{thm}\label{thm:Lpsubloop}
$L_p$ is a subloop of $L$.
\end{thm}
\begin{proof}
For $\c,\d\in L_p$, let $q$ be the greater of the orders of $\c$ and
$\d$, and let $r=q(q-1)/2$.  Then, using di-associativity and the
definition of $\chi$, we have
\begin{equation}
(\c\d)^q = \c^q \d^q \chi(d,c)^r = \chi(d,c)^r,
\end{equation}
and the theorem follows from Theorem~\ref{thm:gcd}.
\end{proof}

\begin{thm}\label{thm:torsionsubloop}
Let $T$ be the set of all elements of $L$ of finite order.  Then $T$
is a subloop of $L$ isomorphic to the restricted direct product of the
$L_p$'s, over all primes $p$.
\end{thm}
\begin{proof}
First, we note that Theorem~\ref{thm:gcd} implies that elements of
relatively prime order commute and associate freely, so the
unassociated product of elements of pairwise relatively prime order is
well-defined.  Consequently, by the Chinese Remainder Theorem, for
every $\c\in L$ of order $n$, we have
\begin{equation}
\c=\prod_{p\dividing n} \c_p,
\end{equation}
where each $\c_p$ is a power of $\c$, and the order of $\c_p$ is a
power of $p$.

It is therefore enough to show that if the orders of $\c_i$ and $\d_j$
are relatively prime for $i,j=1,2$, then $(\c_1\d_1)(\c_2\d_2) =
(\c_1\c_2)(\d_1\d_2)$.  However, using Theorem~\ref{thm:gcd}
repeatedly, we see that 
\begin{equation}
\begin{split}
(\c_1\d_1)(\c_2\d_2)  &= \c_1(\d_1(\c_2\d_2)) = \c_1((\d_1\c_2)\d_2)
	= \c_1((\c_2\d_1)\d_2) \\
	&= \c_1(\c_2(\d_1\d_2)) = (\c_1\c_2)(\d_1\d_2),
\end{split}
\end{equation}
and the theorem follows.
\end{proof}

We next obtain Main Theorem \ref{mainthm:assocsix}.

\begin{proof}[{\bf Proof of Main Theorem \ref{mainthm:assocsix}}]
Since $L^*$ is an abelian group generated by $\alpha(c,d,e)$ for all
$c,d,e\in C$, it is enough to show that $\alpha(c,d,e)^6=1$ for all
$c,d,e\in C$.  However, since \eqref{eq:chimultilin} implies
\begin{equation}
\chi(c,e)\chi(d,e)\alpha(c,d,e)^3 = \chi(cd,e) =
	\chi(dc,e) = \chi(d,e)\chi(c,e)\alpha(d,c,e)^3,
\end{equation}
using skew-symmetry, we have
$\alpha(c,d,e)^3=\alpha(d,c,e)^3=\alpha(c,d,e)^{-3}$, and the theorem
follows.
\end{proof}

We then have the following analogue of Thm.\ 11.2 of Bruck~\cite[Ch.\
VIII]{rhb:loops}.

\begin{thm}\label{thm:finitebyfinitegen}
If $L$ is finitely generated, then $L^*$ is finite.  More precisely,
$L$ is a central extension of $L/L^*$ (a finitely generated group of
class $\le 2$) by a finite group of exponent 6.
\end{thm}
\begin{proof}
Since \eqref{eq:alphamultilin} implies that $L^*$ is an abelian group
generated by $\alpha(c,d,e)$, where $\c$, $\d$, and $\e$ run over all
{\em generators\/} of $L$, the theorem follows from
Main Theorem~\ref{mainthm:assocsix}.
\end{proof}

Finally, we note that Moufang loops of class 2 satisfy the following
stronger version of the Moufang identity.

\begin{thm}\label{thm:LisMklaw}
If $n$ is the exponent of $L^*$, then $L$ satisfies  
\begin{equation}\label{eq:Mklaw}
\c^k(\d(\c\e)) = ((\c^k \d)\c)\e
\end{equation}
for all $\c,\d,\e\in L$, and precisely those integers $k$ such that
$k\=1\pmod{n}$.
\end{thm}
The identity \eqref{eq:Mklaw} is called the {\it $M_k$-law}.  Note
that the $M_1$-law is just \eqref{eq:Mkone}.
\begin{proof}
From the definition of $\alpha$, we have
\begin{equation}
\begin{split}
((\c^k\d)\c)\e &= (\c^k \d)(\c\e) \alpha(c^k d,c,e) \\
	&= \c^k(\d(\c\e)) \alpha(c^k d,c,e)\alpha(c^k,d,ce).
\end{split}
\end{equation}
However, using \eqref{eq:alphasymp}--\eqref{eq:alphamultilin}, we
obtain
\begin{equation}
\begin{split}
\alpha(c^k d,c,e)\alpha(c^k,d,ce)
	&= \alpha(c^k,c,e)\alpha(d,c,e)
		\alpha(c^k,d,c)\alpha(c^k,d,e) \\
	&= \alpha(c,d,e)^{k-1},
\end{split}
\end{equation}
which means that the $M_k$-law is satisfied if and only if the order
of any $\alpha(c,d,e)$ divides $k-1$.  The theorem follows.
\end{proof}

We then have the following corollary.  (This result on Moufang loops
of class 2 can also be obtained more directly from Cor.\ IV.4.8 of
Pflugfelder~\cite{hop:loops}.)

\begin{cor}\label{cor:classtwoGloop}
If $L_3\intersect L^*=1$, then $L$ is a $G$-loop.
\end{cor}
\begin{proof}
If $L_3\intersect L^*=1$, then $L$ satisfies an $M_k$-law for all odd
$k$, and so the corollary follows from Thm.\ IV.4.11 of
Pflugfelder~\cite{hop:loops}.
\end{proof}

\section{Small Frattini Moufang loops and coded vector spaces}
\label{sect:eml}

In the rest of this paper, we assume all loops are finite; in fact, we
will mostly consider loops of prime power order.  To motivate our main
definition (Definition~\ref{defn:SFL}), we begin by imitating Sect.\
23 of Aschbacher~\cite{ma:finite}.

\begin{defn}\label{defn:eml}
Let $p$ be a prime.  We say that a Moufang $p$-loop $L$ is {\it
special\/} if $\Phi(L)=Z(L)=L'$, and we say that a special Moufang
loop $L$ is {\it extraspecial\/} if $Z(L)$ is cyclic.
\end{defn}

For instance, every extraspecial group is an extraspecial Moufang
loop.

Note that every nontrivial special Moufang loop is centrally nilpotent
of class 2.  Theorem~\ref{thm:frattini} therefore implies that if $L$
is a special Moufang loop, then $L/\Phi(L)$ is an elementary abelian
group.  Furthermore, copying the proof of (23.7) in
Aschbacher~\cite{ma:finite} word for word, it also follows that $Z(L)$
is an elementary abelian group.  We conclude that $L$ is an
extraspecial Moufang loop if and only if $\Phi(L)=Z(L)=L'$ has order
$p$.

\begin{rem}\label{rem:whyallp}
As the reader may have noticed, Main Theorem~\ref{mainthm:assocsix}
implies that Definition \ref{defn:eml} is new only when $p=2$ or $3$;
otherwise, we are talking about (extra)special groups.  However, since
it requires little extra effort, we will continue to discuss the case
of arbitrary $p$.
\end{rem}

We generalize our situation slightly with the following definition.

\begin{defn}\label{defn:SFL}
A $p$-loop $L$ is said to be {\it small Frattini\/} if $\Phi(L)$ has
order dividing $p$.  A small Frattini loop $L$ is said to be {\it central
small Frattini\/} if $\Phi(L)\le Z(L)$.
\end{defn}

For instance, every extraspecial Moufang loop is central small
Frattini, as is any elementary abelian group.  More generally:

\begin{thm}\label{thm:SFMLisCSFML}
Every small Frattini Moufang loop is central small Frattini.
\end{thm}
Recall that $C(L)$ denotes the Moufang center of a Moufang loop $L$
(Definition~\ref{defn:Moufangcenter}).
\begin{proof}
Let $L$ be a small Frattini Moufang loop.  The theorem is clear for
groups, so since $L^*\le L'\le \Phi(L)$, we may assume that
$L^*=\Phi(L)$ has order $p$.  It follows that for some $\c,\d,\e\in
L$, $A=[\c,\d,\e]\ne 1$ and $L^*=\gen{A}$.  Now, for all $y\in L$,
$\gen{A}$ is a normal subgroup of order $p$ in $\gen{A,y}$, so
$[A,y]=1$.  In other words, $L^*=\gen{A}\le C(L)$.

Therefore, it is enough to show that $A\in N(L)$.  Furthermore, since
we now know that Proposition \ref{prop:rhblem} applies to $L$, it is
enough to show that $a=[A,\d,\e]=1$ for any $\d,\e\in L$.  However,
for any $\d,\e\in L$, Proposition \ref{prop:rhblem} implies that
$a=[A,\d,\e]$ is central in $\gen{A,\d,\e}$, so if $a\ne 1$, then
$\gen{A}=L^*=\gen{a}=Z(\gen{\c,\d,\e})$, a contradiction.  The theorem
follows.
\end{proof}

\begin{notn}
In the rest of this paper, we abbreviate the term ``small Frattini
Moufang'' as SFM, and we abbreviate ``small Frattini Moufang loop'' as
SFML.  Also, for the rest of this section, let $p$ be a prime, let $L$
be an SFML of order $p^{1+k}$, let $Z$ be a fixed central subgroup of
$L$, and let $C\isom L/Z$ be an elementary abelian $p$-group of rank
$k$ (vector space of dimension $k$ over $\F_p$).  We also retain the
convention of the previous section that $\c,\d,\e,\f\in L$ reduce to
$c,d,e,f\in C$ in the quotient.
\end{notn}

Applying Theorem~\ref{thm:sympmulti}, we see that $\chi$ and $\alpha$
are again well-defined functions which satisfy the formulas
\eqref{eq:chisymp}--\eqref{eq:alphamultilin}.  However, to understand
SFML's, we need one more function.

\begin{defn}\label{defn:sigma}
We define the function $\sigma: C \rightarrow Z$ by
$\sigma(c)=\gamma^p$.  Note that $\sigma$ is well-defined because $Z$
is central and has exponent $p$ and $L/Z$ has exponent $p$.
\end{defn}

The following is the analogue of Theorem \ref{thm:sympmulti} for
$\sigma$.

\begin{thm}\label{thm:linear}
For all $c,d\in C$, we have:
\begin{align}
\sigma(c^n) &= \sigma(c)^n \label{eq:sigmapowerlin} \\
\sigma(cd) &= \begin{cases}
	\sigma(c)\sigma(d)\chi(c,d) &\text{for $p=2$}, \\
	\sigma(c)\sigma(d) &\text{for $p>2$}.
	      \end{cases} \label{eq:sigmalin}
\end{align}
\end{thm}
\begin{proof}
\eqref{eq:sigmapowerlin} is clear.  As for \eqref{eq:sigmalin}, if
$r=p(p-1)/2$, then
\begin{equation}
(\c\d)^p = \c^p \d^p \chi(d,c)^r = \sigma(c)\sigma(d)\chi(d,c)^r.
\end{equation}
The theorem follows from the fact that for $p>2$, $p$ divides $r$, and
for $p=2$, $r=1$ and $\chi(c,d)=\chi(d,c)$.
\end{proof}

We are led to the following definition.

\begin{defn}\label{defn:cvs}
Let $Z$ be the group of order $p$.  A {\it coded vector space\/} (or
CVS) is defined to be a 4-tuple $(C,\sigma,\chi,\alpha)$, where $C$ is
a finite-dimensional vector space over $\F_p$, and $\sigma:
C\rightarrow Z$, $\chi: C\times C\rightarrow Z$, and $\alpha: C\times
C\times C\rightarrow Z$ satisfy
\eqref{eq:sigmapowerlin}--\eqref{eq:sigmalin},
\eqref{eq:chisymp}--\eqref{eq:chimultilin}, and
\eqref{eq:alphasymp}--\eqref{eq:alphamultilin}, for all $c,d,e,f\in C$
and all $n\in\Z$.
\end{defn}

\begin{notn}
We will often refer to the CVS $(C,\sigma,\chi,\alpha)$ simply as
$C$.
\end{notn}

It is worth noting the different forms that \eqref{eq:sigmalin}
and \eqref{eq:chimultilin} take for different $p$.  That is, for
$p=2$, we have
\begin{align}
\sigma(cd) &= \sigma(c)\sigma(d)\chi(c,d),
	\label{eq:sigmalineven} \\
\chi(cd,e) &= \chi(c,e)\chi(d,e)\alpha(c,d,e),
	\label{eq:chimultilineven}
\end{align}
and for $p>2$, we have
\begin{align}
\sigma(cd) &= \sigma(c)\sigma(d),
	\label{eq:sigmalinodd} \\
\chi(cd,e) &= \chi(c,e)\chi(d,e).
	\label{eq:chimultilinodd}
\end{align}
In other words, for $p=2$, $\sigma$, $\chi$, and $\alpha$ are related
by polarization, and for $p>2$, $\sigma$ and $\chi$ are multilinear.
As for $\alpha$, for $p=2$ or 3, $\alpha$ is multilinear, and for
$p>3$, $\alpha$ is identically equal to 1.

We also note that \eqref{eq:sigmapowerlin}, \eqref{eq:chipowerlin},
and \eqref{eq:alphapowerlin} imply
\begin{equation}\label{eq:scaonone}
\sigma(1)=\chi(c,1)=\alpha(c,d,1)=1
\end{equation}
for all $c,d\in C$.

Note that choosing a different generator for $Z$ has the effect of
acting on $\sigma$, $\chi$, and $\alpha$ by an element of $\Aut(Z)$;
in additive terms, this means that $C$ is really only defined up to
scalar multiple.  The natural definition of isomorphism for CVS's is
therefore the following one.

\begin{defn}
Let $(C_i,\sigma_i,\chi_i,\alpha_i)$ be a CVS for $i=1,2$.  We say
that $C_1$ and $C_2$ are {\it isomorphic up to scalar multiple\/} if
there is a vector space isomorphism $\phi: C_1 \rightarrow C_2$ and
some fixed $a\in\Aut(Z)$ such that
\begin{equation}
\begin{split}
\sigma_2(\phi(c)) &= \sigma_1(c)^a, \\
\chi_2(\phi(c),\phi(d)) &= \chi_1(c,d)^a, \\
\alpha_2(\phi(c),\phi(d),\phi(e)) &= \alpha_1(c,d,e)^a.
\end{split}
\end{equation}
If $C_1$ and $C_2$ are isomorphic up to scalar multiple with respect
to the trivial scalar (identity automorphism of $Z$), then we say that
$C_1$ and $C_2$ are {\it isomorphic}.
\end{defn}

Finally, to describe the relationship between SFML's and CVS's, we
introduce one more definition, in which, by convention, we define
$\c^n$ inductively by $\c^0=1$ and $\c^{n+1}=\c\c^n$.

\begin{defn}\label{defn:codedext}
Let $p$ be a prime, and let $(C,\sigma,\chi,\alpha)$ be a CVS over
$\F_p$.  We say that a loop $L$ is a {\it coded extension\/} of $C$ if
$L$ satisfies the following conditions.
\begin{enumerate}
\item\label{cond:centralext} $L$ has a central subgroup $Z$ of order
$p$ such that $L/Z\isom C$.
\item\label{cond:formulas} Let $\c,\d,\e\in L$ denote arbitrary
preimages of $c,d,e\in C$, respectively.  Then:
\begin{align}
\c^p &= \sigma(c),	\label{eq:CEpower} \\
[\c,\d] &= \chi(c,d),	\label{eq:CEcommute} \\
[\c,\d,\e] &= \alpha(c,d,e),	\label{eq:CEassociate}
\end{align}
where the values of $\sigma$, $\chi$, and $\alpha$ are taken to be in
the central subgroup $Z$.
\end{enumerate}
\end{defn}

\begin{thm}\label{thm:codedextisSFML}
Every SFML is a coded extension of a CVS, and every coded extension of
a CVS is an SFML.
\end{thm}
\begin{proof}
If $L$ is a Moufang loop with a central subgroup $Z$ of order $p$ such
that $C\isom L/Z$ is an elementary abelian $p$-group, Theorems
\ref{thm:sympmulti} and \ref{thm:linear} imply that $L$ is a coded
extension of the CVS $(C,\sigma,\chi,\alpha)$, where $\sigma$, $\chi$,
and $\alpha$ are from Definitions \ref{defn:chialpha} and
\ref{defn:sigma}.  Conversely, let $L$ be a coded extension of a CVS
$(C,\sigma,\chi,\alpha)$.  Because $\alpha$ satisfies
\eqref{eq:alphasymp}--\eqref{eq:alphamultilin}, the proof of
Theorem~\ref{thm:LisMklaw} in the case $k=1$ shows that $L$ is
Moufang.  (Note that for $k=1$, the proof of
Theorem~\ref{thm:LisMklaw} does not use power-associativity.)
Therefore, since $L/Z$ is an elementary abelian $p$-group, $L$ is an
SFML.
\end{proof}

\section{Coded vector spaces and doubly even codes}
\label{sect:cvsiscode}

We come to the question: Given an $m$-dimensional vector
space $C$ over $\F_p$, in what ways can $\sigma$, $\chi$, and $\alpha$
be defined to obtain a CVS?  Now, for $p>2$, \eqref{eq:sigmalinodd}
and \eqref{eq:chimultilinodd} show that $\sigma$, $\chi$, and $\alpha$
can be chosen independently, which makes this question easy.
On the other hand, for $p=2$, $\sigma$, $\chi$, and $\alpha$ are
related by polarization, so it is less clear {\it a priori\/} which
CVS's exist over $\F_2$.

Let $C$ be a CVS over $\F_2$, and let $\set{c_1\dots c_m}$ be a basis
for $C$.  Clearly, the symplectic, skew-symmetric (or in
characteristic 2, symmetric), and polarization properties of $\sigma$,
$\chi$, and $\alpha$ imply that $\sigma$, $\chi$, and $\alpha$ are
determined by $\sigma(c_i)$ ($1\le i\le m$), $\chi(c_i,c_j)$ ($1\le
i<j \le m$), and $\alpha(c_i,c_j,c_k)$ ($1\le i<j<k \le m$).
Conversely, as we shall see in a moment, we may define a valid
$\sigma$, $\chi$, and $\alpha$ by setting these values arbitrarily.
Now, it is possible to prove this directly (see Theorem
\ref{thm:consistenttwo}), but we will instead show that every CVS over
$\F_2$ can be obtained from a {\it doubly even code\/}
(Theorem~\ref{thm:CVSiscode}).

\begin{notn}
For the rest of this section, we revert to additive notation for
$\F_2$.
\end{notn}

\begin{defn}\label{defn:code}
A {\it binary code\/} (or in this paper, simply a code) of length $n$
and dimension $m$ is defined to be a subspace of $\F_2^n$ of dimension
$m$.  We define $\wt{c}$ (resp.\ $\wt{c\intersect d}$,
$\wt{c\intersect d\intersect e}$) to be the number of non-0
coordinates in $c$ (resp.\ common to $c$ and $d$, common to $c$, $d$,
and $e$).  We say that a code $C$ is {\it doubly even\/} if
$\wt{c}\=0$ (mod~4) for all $c\in C$.
\end{defn}

Note that if $C$ is doubly even, then $\wt{c\intersect d}\=0\pmod{2}$
for all $c,d\in C$.  Conversely, if $\set{c_i}$ is a basis for a code
$C$, it is easy to see that $C$ is doubly even if and only if
$\wt{c_i}\=0\pmod{4}$ and $\wt{c_i\intersect c_j}\=0\pmod{2}$ for all
$c_i$ and $c_j$ in the basis.

Doubly even codes determine CVS's in the following manner.

\begin{thm}\label{thm:codeisCVS}
Let $C$ be a doubly even code, and define
\begin{align}
\sigma(c) &\= \frac{\wt{c}}{4} \pmod{2}, \\
\chi(c,d) &\= \frac{\wt{c\intersect d}}{2} \pmod{2}, \\
\alpha(c,d,e) &\= \wt{c\intersect d\intersect e} \pmod{2},
\end{align}
for all $c,d,e\in C$.  Then $(C,\sigma,\chi,\alpha)$ is a CVS.
\end{thm}
\begin{proof}
The symmetry of $\chi$ and $\alpha$ is clear, as are
\eqref{eq:sigmapowerlin}, \eqref{eq:chipowerlin}, and
\eqref{eq:alphapowerlin}.  Furthermore, $\chi$, resp.\ $\alpha$, is
symplectic (\eqref{eq:chisymp}, resp.\ \eqref{eq:alphasymp}) because
$\wt{c}\=0\pmod{4}$, resp.\ $\wt{c\intersect d}\=0\pmod{2}$, for all
$c,d\in C$.  As for polarization (\eqref{eq:sigmalineven},
\eqref{eq:chimultilineven}, and \eqref{eq:alphamultilin}), suppose we
associate with each $c\in C$ a diagonal matrix $M_c$ with integer
entries of 0's and 1's corresponding to the coordinates of $c$.  Then
$\wt{c}=\tr M_c$, $\wt{c\intersect d}=\tr \(M_c M_d\)$, and
$\wt{c\intersect d\intersect e}=\tr \(M_c M_d M_e\)$, which means that
\begin{equation}
M_{cd} = M_c + M_d - 2 M_c M_d
\end{equation}
implies polarization.
\end{proof}

Conversely, we have:

\begin{thm}\label{thm:CVSiscode}
For any integer $m>0$, choose elements $\sigma_i$ ($1\le i\le m$),
$\chi_{ij}$ ($1\le i<j \le m$), and $\alpha_{ijk}$ ($1\le i<j<k \le
m$) of $\F_2$.  There exists a (unique) CVS $C=\gen{c_1\dots c_m}$ of
dimension $m$ such that $\sigma(c_i)=\sigma_i$ ($1\le i\le m$),
$\chi(c_i,c_j)=\chi_{ij}$ ($1\le i<j \le m$), and
$\alpha(c_i,c_j,c_k)=\alpha_{ijk}$ ($1\le i<j<k \le m$).  Furthermore,
$C$ is isomorphic (as a CVS) to a doubly even code.
\end{thm}
\begin{proof}
Proceeding by induction on $m$, it suffices to construct a code $C$
which has the required values of $\sigma$, $\chi$, and $\alpha$ on its
basis.  For $m=1$, if $\sigma_1=0$, we may take
$C=\gen{(1,1,1,1,1,1,1,1)}$, and if $\sigma_1=1$, we may take
$C=\gen{(1,1,1,1)}$.  By induction, then, let $C_0=\gen{c_1\dots
c_{m-1}}$ be a doubly even code with the correct values of $\sigma$,
$\chi$, and $\alpha$ on all basis vectors with indices $<m$.  We
construct $C$ using the following steps.
\begin{enumerate}
\item If $\sigma_m=0$ (resp.\ $\sigma_m=1$), extend $C_0$ by 8 (resp.\
4) coordinates, extend the basis vectors $c_1\dots c_{m-1}$ by 0's,
and add a new vector $c_m$ with all coordinates 0 except the last 8
(resp.\ 4).  The resulting doubly even code $C_1$ then has the correct
values of $\sigma$ on the basis, and the correct values for $\chi$ and
$\alpha$ on the basis when all of the indices are $<m$; however,
$\chi(c_i,c_m)=0$ for $1\le i<m$ and $\alpha(c_i,c_j,c_m)=0$ for $1\le
i<j<m$.
\item For each $i<m$ such that $\chi_{im}=1$, extend $C_1$ by 14
coordinates, extend $c_i$ and $c_m$ by
\begin{equation}
\begin{array}{rcccccccccccccc}
c_i\quad
	1 & 1 & 1 & 1 & 1 & 1 & 1 & 1 & 0 & 0 & 0 & 0 & 0 & 0 \\
c_m\quad
	0 & 0 & 0 & 0 & 0 & 0 & 1 & 1 & 1 & 1 & 1 & 1 & 1 & 1,
\end{array}
\end{equation}
and extend the other basis vectors by 0's.  The resulting code $C_2$
has the correct values of $\sigma$ and $\chi$ on the basis, so it
remains to correct $\alpha$.
\item For each $i<j<m$ such that $\alpha_{ijm}=1$, extend $C_2$ by 13
coordinates, extend $c_i$, $c_j$, and $c_m$ by
\begin{equation}
\begin{array}{rccccccccccccc}
c_i\quad
	1 & 1 & 1 & 1 & 0 & 0 & 0 & 1 & 1 & 1 & 1 & 0 & 0 \\
c_j\quad
	1 & 1 & 1 & 1 & 1 & 1 & 1 & 0 & 0 & 0 & 0 & 1 & 0 \\
c_m\quad
	1 & 0 & 0 & 0 & 1 & 1 & 1 & 1 & 1 & 1 & 0 & 0 & 1,
\end{array}
\end{equation}
and extend the other basis vectors by 0's.  The resulting code $C$
has the correct values of $\sigma$, $\chi$, and $\alpha$ on the basis,
and the theorem follows.  \qed
\end{enumerate}
\renewcommand{\qed}{}
\end{proof}

Main Theorem \ref{mainthm:evenSFMLiscodeloop} now follows.

\begin{proof}[{\bf Proof of Main Theorem
	\ref{mainthm:evenSFMLiscodeloop}}]
If $C$ is a doubly even code, it is easy to see that the coded
extensions of $C$ are precisely the code loops over $C$.  (Compare
Conway~\cite{jhc:monster} and Griess~\cite{rlg:codeloops}.)
Therefore, Main Theorem \ref{mainthm:evenSFMLiscodeloop} follows from
Theorems \ref{thm:codedextisSFML}, \ref{thm:codeisCVS}, and
\ref{thm:CVSiscode}.
\end{proof}

\begin{rem}\label{rem:griesspf}
Griess (personal communication) has suggested the following alternate
proof for Main Theorem \ref{mainthm:evenSFMLiscodeloop}.  Now, it is
not hard to see that the function $\sigma$ of a CVS over $\F_2$ is a
Parker function, as described in Defn.\ 13 of
Griess~\cite{rlg:codeloops}.  (For instance, use \eqref{eq:sigmatwo},
below.)  This fact and Theorem \ref{thm:codedextisSFML} imply that
every SFM 2-loop is a loop constructed by ``Parker's procedure'',
which means that Main Theorem \ref{mainthm:evenSFMLiscodeloop} follows
from Thm.\ 14 of Griess~\cite{rlg:codeloops}.  In other words, in this
approach, Griess' Thm.\ 14 replaces Theorems
\ref{thm:codeisCVS} and \ref{thm:CVSiscode}.
\end{rem}

\begin{rem}\label{rem:codesize}
Note that there are many ways of expressing a given SFM 2-loop as a
code loop.  For instance, let $\set{c,d,e}$ be a basis for the CVS of
the octonion loop (the standard double basis for the octonions).
Then, in additive notation,
\begin{equation}
\sigma(c)=\sigma(d)=\sigma(e)=\chi(c,d)=\chi(c,e)=\chi(d,e)
	=\alpha(c,d,e)=1.
\end{equation}
Therefore, applying the procedure in the proof of Theorem
\ref{thm:CVSiscode}, we obtain a doubly even code of length
$3(4)+3(14)+13=67$, the code loop of which is the octonions.  On the
other hand, applying the procedure in Thm.\ 14 of
Griess~\cite{rlg:codeloops} to the Parker loop structure of the
octonions gives a code of length 38.  It is also easily shown that the
octonions are the code loop over the Hamming $[7,3,4]$ code.
\end{rem}

\section{Existence and uniqueness of coded extensions}
\label{sect:codedextensions}

Because of Theorem \ref{thm:codedextisSFML}, to obtain Main Theorem
\ref{mainthm:CVSisSFML}, it remains to show that every CVS has a
unique coded extension, which we do in this section (Theorems
\ref{thm:codedextunique} and \ref{thm:codedextexist}).  Now, for code
loops ($p=2$), this is Thm.\ 10 of Griess~\cite{rlg:codeloops}, and
for $p>3$, we are in the associative case, which means that this is
essentially known.  However, to gain insight into the known cases and
to introduce the {\it semidirect central product\/}
(Definition~\ref{defn:sdcp}), we continue to consider all primes $p$.

\begin{notn}
In this section, we use script letters $\C,\D,\E,\dots$ to denote
coded extensions of the CVS's $C,D,E,\dots$, possibly with subscripts.
\end{notn}

We first address uniqueness.  Our proof follows \S2 of
Conway~\cite{jhc:monster}.

\begin{thm}\label{thm:codedextunique}
For $n=1,2$, let $C_n$ be a CVS of dimension $k$ over $\F_p$, and let
$\C_n$ be a coded extension of $C_n$.  If $C_1$ and $C_2$ are
isomorphic up to scalar multiple, then $\C_1$ is isomorphic to $\C_2$.
\end{thm}
\begin{proof}
First, by choosing a different generator for the distinguished central
subgroup of $\C_2$, we may assume that $C_1$ and $C_2$ are isomorphic.
So let $\set{c_i}$ be a basis for $(C,\sigma,\chi,\alpha)=C_1$, and
let $\C$ be the loop given by the loop presentation
\begin{align}
\left\langle\vphantom{\beta^{-1}}\right. z, \gamma_i
        	\left.\vphantom{\beta^{-1}}\right|
		\qquad\qquad& \notag \\
	\c_i^p &= \sigma(c_i),
		\label{eq:Cpowerreln} \\
	[\c,\d] &= \chi(c,d),
		\label{eq:Ccommreln} \\
	[\c,\d,\e] &= \alpha(c,d,e),
		\label{eq:Cassocreln} \\
        z^p &= 1
		\label{eq:zpower}
                \left.\vphantom{\beta^{-1}}\right\rangle,
\end{align}
where $i$ runs between $1$ and $k$; $\c,\d,\e$ run over all loop words
in the generators; $c,d,e$ are the images in $C$ of $\c,\d,\e$,
respectively, under the map sending $\c_i$ to $c_i$ and $z$ to 1; the
values of $\sigma$, $\chi$, and $\alpha$ are taken to be in
$Z=\gen{z}$, which is a central subgroup of $\C$ because of
\eqref{eq:Ccommreln}--\eqref{eq:zpower} and \eqref{eq:scaonone}; and
the expression $\c_i^p$ is defined inductively by $\c_i^0=1$ and
$\c_i^{n+1}=\c_i\c_i^n$.

Now, since isomorphisms take bases to bases, the above presentation is
purely a function of the isomorphism class of $C$.  Therefore, it is
enough to show that any coded extension of $C$ is isomorphic to $\C$.
Furthermore, since the universal property of loop presentations (see
Evans~\cite[I.2]{te:varloops}) implies that {\em any\/} coded
extension of $C$ is a homomorphic image of $\C$, it is enough to show
that $\C$ has order at most $p^{1+k}$.  However, since
\eqref{eq:Ccommreln}--\eqref{eq:zpower} imply that $Z$ is a central
subgroup of $\C$ of order dividing $p$, and
\eqref{eq:Cpowerreln}--\eqref{eq:Cassocreln} implies that $\C/Z$ is an
elementary abelian $p$-group of rank $k$, the theorem follows.
\end{proof}

\begin{notn}
We resume the convention that if $\c,\d,\e$ (possibly with subscripts)
are elements of a coded extension, then $c,d,e$ (possibly with
subscripts) are the corresponding elements of the quotient CVS.
\end{notn}

We turn to existence.  Now, if $\D$ and $\E$ are subloops of the same
coded extension, then for $\d_i\in\D$, $\e_i\in\E$,
$(\d_1\e_1)(\d_2\e_2)=z_0(\d_1\d_2)(\e_1\e_2)$, where $z_0\in Z$ is
expressible in terms of $\chi$, $\alpha$, $d_1$, $e_1$, $d_2$ and
$e_2$.  This observation motivates the following defintion.

\begin{defn}\label{defn:sdcp}
Let $p$ be a prime and let $Z=\F_p$ be the group of order $p$.  Let
$D$ and $E$ be linearly independent subspaces of a CVS
$(C,\sigma,\chi,\alpha)$; note that $D$ and $E$ are CVS's by
restriction.  Let $\D$ (resp.\ $\E$) be a coded extension of $D$
(resp.\ $E$).  It is easily verified that the binary operation on the
set $Z\times\D\times\E$ given by
\begin{equation}\label{eq:semicentral}
(z_1,\d_1,\e_1)(z_2,\d_2,\e_2) = (z_0 z_1 z_2,\d_1\d_2,\e_1\e_2),
\end{equation}
where
\begin{equation}\label{eq:sdcperror}
z_0 = \chi(e_1,d_2) \alpha(d_1,e_1 d_2\i,e_2)
	\alpha(d_1,e_1,d_2)^2 \alpha(e_1,d_2,e_2)^{-2},
\end{equation}
defines a loop $\Gamma$ containing the central subgroup $Z\times
Z\times Z$.  We may therefore define the {\it semidirect central
product\/} of $\D$ and $\E$ to be the quotient of $\Gamma$ by the
central subgroup $\gen{(z,z\i,1),(z,1,z\i)}$, where $z$ runs over all
elements of $Z$.
\end{defn}

Note that for $p=2$, \eqref{eq:sdcperror}
becomes
\begin{equation}\label{eq:semicentraleven}
z_0 = \chi(e_1,d_2) \alpha(d_1,e_1 d_2,e_2),
\end{equation}
for $p=3$, \eqref{eq:sdcperror} becomes
\begin{equation}\label{eq:semicentralodd}
z_0 = \chi(e_1,d_2) \alpha(d_1,e_1 d_2\i,e_2)
	\alpha(d_1,e_1,d_2)\i \alpha(e_1,d_2,e_2),
\end{equation}
and for $p>3$,
\begin{equation}\label{eq:semicentralbigp}
z_0 = \chi(e_1,d_2).
\end{equation}

\begin{notn}
Imitating the ATLAS~\cite{ATLAS} notation $\D\cp\E$ for the central
product of two groups $\D$ and $\E$, we use $\D\sdcp\E$ to denote the
semidirect central product (``central product up to sign'') of two
code loops $\D$ and $\E$.  Note that if the $\chi$ and $\alpha$
factors on the right-hand side of \eqref{eq:semicentral} are always
$1$, then $\D\sdcp\E$ just becomes the ordinary central product.  In
this situation, we write $\D\cp\E$, just as in the group case.
\end{notn}

\begin{thm}\label{thm:sdcpiscodedext}
Let $D$ and $E$ be linearly independent sub-CVS's of
$(C,\sigma,\chi,\alpha)$, let $\D$ (resp.\ $\E$) be a coded extension
of $D$ (resp.\ $E$), and let $\C=\D\sdcp\E$.  Then $\C$ is a coded
extension of $D\oplus E$.
\end{thm}
Throughout the following proof, by convention, $\c=(z,\d,\e)$,
possibly with subscripts.  (Note that by our usual convention, we then
have $c=de$.)  We also freely identify $Z$ with $Z\times 1\times 1$.
Finally, the ``skew-symmetric,'' ``symplectic,'' and ``multilinear''
properties of $\chi$ and $\alpha$
(\eqref{eq:chisymp}--\eqref{eq:alphamultilin}) will be applied freely.
\begin{proof}
First, taking the quotient of $Z\times Z\times Z$ in $\C$ as our
distinguished central subgroup, it is easy to see that condition
\ref{cond:centralext} of Definition~\ref{defn:codedext} holds in $\C$.
Furthermore, by collecting ``signs'' (elements of $Z$), it is easy to
see that \eqref{eq:CEpower}--\eqref{eq:CEassociate} hold ``up to
sign,'' so it remains to check the signs.

We first verify \eqref{eq:CEpower}, recalling our convention that
$\c^n$ is defined inductively by $\c^0=1$ and $\c^{n+1}=\c\c^n$.
First, we claim that for all $n\ge 0$,
\begin{equation}\label{eq:powerind}
\c^n = (z^n\chi(e,d)^r,\d^n,\e^n)
\end{equation}
where $r=n(n-1)/2$.  In fact, if \eqref{eq:powerind} holds, then
\begin{equation}
\begin{split}
\c^{n+1} &= \c\c^n \\
	&= (z,\d,\e)(z^n\chi(e,d)^r,\d^n,\e^n) \\
	&= (z^{n+1}\chi(e,d)^r\chi(e,d^n),\d^{n+1},\e^{n+1}) \\
	&= (z^{n+1}\chi(e,d)^{r+n},\d^{n+1},\e^{n+1}),
\end{split}
\end{equation}
and since $r+n=n(n+1)/2$, \eqref{eq:powerind} follows by induction.
In particular, if $n=p$, $r=p(p-1)/2$, which means that
$\chi(e,d)^r=\chi(d,e)$ for $p=2$ and $\chi(e,d)^r=1$ for $p>2$.
Therefore, collecting signs, we get
\begin{equation}
\begin{split}
\c^p &= (\chi(e,d)^r,\sigma(d),\sigma(e)) \\
	&= \begin{cases}
	\sigma(d)\sigma(e)\chi(d,e) &\text{for $p=2$} \\
	\sigma(d)\sigma(e) &\text{for $p>2$}
	   \end{cases} \\
	&= \sigma(de) \\
	&= \sigma(c).
\end{split}
\end{equation}

Next, turning to \eqref{eq:CEcommute}, by collecting signs, we see
that
\begin{equation}
\c_1\c_2 = (z_3 z_1 z_2,\d_1\d_2,\e_1\e_2),
\end{equation}
where
\begin{equation}
z_3 = \chi(e_1,d_2) \alpha(d_1,e_1 d_2\i,e_2) \alpha(d_1,e_1,d_2)^2
	\alpha(e_1,d_2,e_2)^{-2};
\end{equation}
and 
\begin{equation}
\c_2\c_1 = (z_4 z_1 z_2,\d_1\d_2,\e_1\e_2),
\end{equation}
where
\begin{equation}\label{eq:zfour}
\begin{split}
z_4 &= \chi(e_2,d_1) \chi(d_2,d_1) \chi(e_2,e_1) \\
	&\qquad \alpha(d_2,e_2 d_1\i,e_1) \alpha(d_2,e_2,d_1)^2
	\alpha(e_2,d_1,e_1)^{-2}.
\end{split}
\end{equation}
To obtain \eqref{eq:CEcommute}, we need to show that $z = z_3 z_4\i =
\chi(d_1 e_1,d_2 e_2)$.  However, gathering the $\chi$ terms of $z$,
we have
\begin{equation}
\begin{split}
&\chi(e_1,d_2)\chi(e_2,d_1)\i\chi(d_2,d_1)\i\chi(e_2,e_1)\i \\
&\qquad = \chi(e_1,d_2)\chi(d_1,e_2)\chi(d_1,d_2)\chi(e_1,e_2) \\
&\qquad = \chi(d_1,d_2 e_2)\chi(e_1,d_2 e_2)
	\alpha(d_1,d_2,e_2)^3 \alpha(e_1,d_2,e_2)^3 \\
&\qquad = \chi(d_1 e_1,d_2 e_2)
	\alpha(d_1,e_1,d_2 e_2)^3 \alpha(d_1 e_1,d_2,e_2)^3,
\end{split}
\end{equation}
and gathering the $\alpha$ terms of $z$, we have
\begin{equation}
\begin{split}
&
	\alpha(d_1,e_1 d_2\i,e_2)
	\alpha(d_1,e_1,d_2)^2\alpha(e_1,d_2,e_2)^{-2} \\
{}\cdot{} &
	\alpha(d_2,e_2 d_1\i,e_1)\i
	\alpha(d_2,e_2,d_1)^{-2}\alpha(e_2,d_1,e_1)^2 \\
{}={} &
	\alpha(d_1,e_1,d_2)^3\alpha(d_1,e_1,e_2)^3
	\alpha(d_1,d_2,e_2)^3\alpha(e_1,d_2,e_2)^3.
\end{split}
\end{equation}
\eqref{eq:CEcommute} follows because $\alpha^6=1$ identically.

Finally, using the same strategy to verify \eqref{eq:CEassociate}, we
see that $(\c_1\c_2)\c_3 = z \cdot \c_1(\c_2\c_3)$, where
\begin{equation*}
\begin{split}
z {}={} &
	\chi(e_1,d_2) \alpha(d_1,e_1 d_2\i,e_2)
	\alpha(d_1,e_1,d_2)^2 \alpha(e_1,d_2,e_2)^{-2} \\
{}\cdot{} &
	\chi(e_1 e_2,d_3) \alpha(d_1 d_2,e_1 e_2 d_3\i,e_3)
	\alpha(d_1 d_2,e_1 e_2,d_3)^2 \alpha(e_1 e_2,d_3,e_3)^{-2} \\
{}\cdot{} &
	\chi(e_2,d_3)\i \alpha(d_2,e_2 d_3\i,e_3)\i
	\alpha(d_2,e_2,d_3)^{-2} \alpha(e_2,d_3,e_3)^2 \\
{}\cdot{} &
	\chi(e_1,d_2 d_3)\i \alpha(d_1,e_1 d_2\i d_3\i,e_2 e_3)\i
	\alpha(d_1,e_1,d_2 d_3)^{-2} \alpha(e_1,d_2 d_3,e_2 e_3)^2 \\
{}\cdot{} &
	\alpha(d_1,d_2,d_3)\alpha(e_1,e_2,e_3).
\end{split}
\end{equation*}

Collecting the $\chi$ terms of $z$, we have
\begin{equation}\label{eq:chitermsofalpha}
\begin{split}
&\chi(e_1,d_2)\chi(e_1 e_2,d_3)\chi(e_2,d_3)\i\chi(e_1,d_2 d_3)\i \\
&\qquad
	\begin{split}
	{}={}
	&\chi(e_1,d_2)\chi(e_1,d_3)\chi(e_2,d_3)
		\alpha(e_1,e_2,d_3)^3 \\
	&\chi(e_2,d_3)\i\chi(e_1,d_2)\i\chi(e_1,d_3)\i
		\alpha(e_1,d_2,d_3)^3
	\end{split} \\
&\qquad = \alpha(e_1,e_2,d_3)^3 \alpha(e_1,d_2,d_3)^3.
\end{split}
\end{equation}

Completely expanding all $\alpha$ terms of $z$ and collecting all
terms which repeat subscripts, we have
\begin{equation}\label{eq:magiccancel}
\begin{split}
%% (12)3 terms
&
	\alpha(d_1,e_1,e_2)\alpha(d_1,d_2\i,e_2)
	\alpha(d_1,e_1,d_2)^2\alpha(e_1,d_2,e_2)^{-2} \\
{}\cdot{} &
	\alpha(d_1,d_3\i,e_3)\alpha(d_2,d_3\i,e_3)
	\alpha(d_1,e_1,e_3)\alpha(d_2,e_2,e_3) \\
{}\cdot{} &
	\alpha(d_1,e_1,d_3)^2\alpha(d_2,e_2,d_3)^2
	\alpha(e_1,d_3,e_3)^{-2}\alpha(e_2,d_3,e_3)^{-2} \\[1ex]
%% space has been inserted to separate two cancelling halves of
%% equation
%% 1(23) terms
{}\cdot{} &
	\alpha(d_2,e_2,e_3)\i\alpha(d_2,d_3\i,e_3)\i
	\alpha(d_2,e_2,d_3)^{-2}\alpha(e_2,d_3,e_3)^2 \\
{}\cdot{} &
	\alpha(d_1,e_1,e_2)\i\alpha(d_1,e_1,e_3)\i
	\alpha(d_1,d_2\i,e_2)\i\alpha(d_1,d_3\i,e_3)\i \\
{}\cdot{} &
	\alpha(d_1,e_1,d_2)^{-2}\alpha(d_1,e_1,d_3)^{-2}
	\alpha(e_1,d_2,e_2)^2\alpha(e_1,d_3,e_3)^2
	= 1,
\end{split}
\end{equation}
since (thankfully) every term in the first three lines of
\eqref{eq:magiccancel} is cancelled by some term in the next three
lines.

The remaining $\alpha$ terms of $z$ are then, after expansion,
\begin{equation}\label{eq:allsubs}
\begin{split}
&
	\alpha(d_1,e_2,e_3)\alpha(d_2,e_1,e_3)
	\alpha(d_1,e_2,d_3)^2\alpha(d_2,e_1,d_3)^2 \\
{}\cdot{} &
	\alpha(d_1,d_2\i,e_3)\i\alpha(d_1,d_3\i,e_2)\i
	\alpha(e_1,d_2,e_3)^2\alpha(e_1,d_3,e_2)^2 \\
{}\cdot{} &
	\alpha(d_1,d_2,d_3)\alpha(e_1,e_2,e_3) \\
{}={} &
	\alpha(d_1,e_2,e_3)\alpha(e_1,d_2,e_3)\i
	\alpha(d_1,e_2,d_3)^2\alpha(e_1,d_2,d_3)^{-2} \\
{}\cdot{} &
	\alpha(d_1,d_2,e_3)\alpha(d_1,e_2,d_3)\i
	\alpha(e_1,d_2,e_3)^2\alpha(e_1,e_2,d_3)^{-2} \\
{}\cdot{} &
	\alpha(d_1,d_2,d_3)\alpha(e_1,e_2,e_3) \\
{}={} &
	\alpha(d_1,e_2,e_3)\alpha(d_1,e_2,d_3)
	\alpha(e_1,d_2,d_3)^{-2} \\
{}\cdot{} &
	\alpha(d_1,d_2,e_3)\alpha(e_1,d_2,e_3)
	\alpha(e_1,e_2,d_3)^{-2} \\
{}\cdot{} &
	\alpha(d_1,d_2,d_3)\alpha(e_1,e_2,e_3).
\end{split}
\end{equation}
Multiplying \eqref{eq:chitermsofalpha}--\eqref{eq:allsubs}, we see
that $z=\alpha(d_1 e_1,d_2 e_2,d_3 e_3)$, so \eqref{eq:CEassociate},
and the theorem, follow.
\end{proof}

\begin{rem}
Compare Kitazume~\cite{mk:codeloops}, whose Thm.\ 2 is a particular
example of the $p=2$ case of the above theorem.
\end{rem}

With Theorem \ref{thm:sdcpiscodedext}, to obtain the existence of
coded extensions, we just need the following example.

\begin{exmp}\label{exmp:dimone}
Let $(C,\sigma,1,1)$ be a CVS of dimension 1, and let $c$ be a nonzero
vector in $C$.  If $\sigma(c)=1$, then $Z\times C$ is a coded
extension of $C$; otherwise, the cyclic group of order $p^2$ is a
coded extension of $C$.
\end{exmp}

\begin{thm}\label{thm:codedextexist}
If $C$ is a CVS, there is a (unique) coded extension of $C$.
\end{thm}
\begin{proof}
Proceeding by induction on $k=\dim C$, let $C=D+E$, where $\dim D=k-1$
and $\dim E=1$, let $\D$ be the coded extension of $D$ (by induction),
and let $\E$ be the coded extension of $E$ (from Example
\ref{exmp:dimone}).  Then from Theorem \ref{thm:sdcpiscodedext},
$\D\sdcp\E$ is a coded extension of $C$.
\end{proof}

\begin{rem}\label{rem:sdcp}
Consider again the loop $\C$ given by the presentation in the proof of
Theorem \ref{thm:codedextunique}.  Now, there is an easy ``solution''
to the loop word problem for $\C$.  Namely, given a loop word $w$ in
the generators $\set{\c_i}$, $w$ can be arranged into the normal form
$z\c_1^r(\c_2^s(\c_3^t(\dots)))$ by simply powering, commuting, and
associating elements, while keeping track of ``error terms'' with
$\sigma$, $\chi$, and $\alpha$.  Theorem \ref{thm:codedextexist} can
then be interpreted as saying precisely that this solution is actually
consistent, i.e., that we are never forced to kill $Z$.

In fact, this normal form procedure provides a method for doing
calculations in any SFML.  Furthermore, while building an SFML from
1-dimensional pieces becomes unwieldy for large dimension, by using
larger associative (or at least familiar) pieces as building blocks,
and gluing them together with the semidirect central product, it is
not hard to do hand calculations in, say, 12- or 16-dimensional
examples.  See~\cite{th:codeloops}.
\end{rem}

\begin{rem}\label{rem:conway}
Conway (personal communication) has remarked that he and Par\-ker
orignally verified the existence of the Parker loop (unpublished) by
proving a ``consistency'' theorem analogous to MacLane's theorems for
monoidal categories~\cite[Ch.\ VII]{sm:categories}, using
``pentagonal'' and ``hexagonal'' relations in the Parker loop.  This
approach could also probably be used to prove Theorem
\ref{thm:codedextexist}.
\end{rem}

Finally, we prove Main Theorem \ref{mainthm:CVSisSFML}.

\begin{proof}[{\bf Proof of Main Theorem \ref{mainthm:CVSisSFML}}]
Because of Theorems \ref{thm:codedextunique} and
\ref{thm:codedextexist}, it remains only to show that if $L$ and $M$
are SFML's, with distinguished central subgroups $Z_L$ and $Z_M$, and
associated CVS's $C_L$ and $C_M$, and there is an isomorphism $\phi: L
\rightarrow M$ such that $\phi(Z_L)=Z_M$, then $C_L$ and $C_M$ are
isomorphic up to scalar multiple.  However, since $\sigma$, $\chi$,
and $\alpha$ are only determined by the isomorphism type of an SFML,
by choosing corresponding generators of $Z_L$ and $Z_M$ (i.e., by
applying a scalar multiple), we can make $C_L$ and $C_M$ isomorphic.
The theorem follows.
\end{proof}

\section{Isotopy in small Frattini Moufang 3-loops}
\label{sect:isotopy}

In this section, as a fairly straightforward application of Main
Theorem \ref{mainthm:CVSisSFML}, we characterize isotopy in SFML's.
Now, because of Corollary \ref{cor:classtwoGloop}, any SFM $p$-loop is
a $G$-loop, unless $p=3$.  Therefore, even though much of what we say
will apply in general, throughout this section, we assume that $L$ is
a finite SFM 3-loop, that $Z$ is a fixed central subgroup of $L$, and
that $C\isom L/Z$ is an elementary abelian 3-group.  In other words,
we assume that $L$ is a coded extension of a CVS
$(C,\sigma,\chi,\alpha)$ over $\F_3$.  We also resume the convention
that $\c,\d,\e,\kappa\in L$ reduce to $c,d,e,k\in C$.

\begin{defn}
Let $L$ be a Moufang loop.  For $\c,\d,\kappa\in L$, we define
\begin{equation}
\c \circ_\kappa \d = (\c\kappa)(\kappa\i\d).
\end{equation}
It is easily verified that the set $L$ and the operation
$\circ_\kappa$ form a loop $(L,\circ_\kappa)$, with identity 1.
Because $(U,V,\iota)$ is an isotopism from $L$ to $(L,\circ_\kappa)$,
where $U(x)=x\kappa\i$, $V(x)=\kappa x$, and $\iota(x)=x$,
$(L,\circ_\kappa)$ is called the {\it $\kappa$-isotope\/} of $L$.
\end{defn}

\begin{rem}\label{rem:circinv}
Note that the di-associativity of $L$ implies that for $\c\in L$,
$\c\i$ is still the inverse of $\c$ in $(L,\circ_\kappa)$.
\end{rem}

We can now state the following fundamental result on isotopy in
Moufang loops (see Pflugfelder~\cite[Thm.\ IV.4.1]{hop:loops}).

\begin{thm}\label{thm:isoiskappa}
If $L$ is a Moufang loop, then any loop-isotope of $L$ is isomorphic
to a $\kappa$-isotope of $L$.  \qed
\end{thm}

The following formula therefore reduces isotopy in SFML's to a matter
of multilinear algebra.

\begin{thm}\label{thm:kappaformula}
We have
\begin{equation}\label{eq:kappa}
\c \circ_\kappa \d = \c\d\cdot\alpha(c,k,d).
\end{equation}
\end{thm}
\begin{proof}
Applying the multilinear and symplectic properties of $\alpha$, we
have:
\begin{equation}
\begin{split}
\c \circ_\kappa \d &= (\c\kappa)(\kappa\i\d) \\
	&= \c(\kappa(\kappa\i\d)) \cdot\alpha(c,k,k\i d) \\
	&= \c\d \cdot\alpha(c,k,d).  \qed
\end{split}
\end{equation}
\noqed\end{proof}

We now only need the following definitions to proceed.

\begin{defn}
The {\it radical\/} of $\alpha$ (resp.\ $\chi$), denoted by
$\rad(\alpha)$ (resp.\ $\rad(\chi)$), is the set of all $c\in C$ such
that $\alpha(c,d,e)=1$ (resp.\ $\chi(c,d)=1$) for all $d,e\in C$
(resp.\ $d\in C$).  Note that $\rad(\alpha) = N(L)/Z$ and $\rad(\chi)
= C(L)/Z$.
\end{defn}

\begin{defn}
Let $(C,\sigma(c),\chi(c,d),\alpha(c,d,e))$ be a CVS.  For $k\in C$,
the {\it adjoint translate\/} $\adt_k(C)$ is defined to be the CVS
$(C,\sigma(c),\chi(c,d)\alpha(c,k,d),\alpha(c,d,e))$.  The operation
$\adt_k$ is called an {\it adjoint translation}.
\end{defn}

%%\begin{rem}
%%Note that for $k,m\in C$, multlinearity of $\alpha$ implies that
%%\begin{equation}
%%\adt_{km}(C) = \adt_k(\adt_m(C)).
%%\end{equation}
%%In other words, the adjoint translations on a CVS
%%$(C,\sigma,\chi,\alpha)$ form a group naturally isomorphic to
%%$C/\rad(\alpha)$.
%%\end{rem}

In the following, we write $\c\circ_\kappa\d$ as $\c\circ\d$.  Also,
for $\c,\d,\e\in L$, we define $\c^{\circ n}$ inductively by
$\c^{\circ 0}=1$ and $\c^{\circ (n+1)} = \c\circ\c^{\circ n}$; we
define $[\c,\d]_\circ$ to be $(\d\circ\c)\i(\c\circ\d)$; and we define
$[\c,\d,\e]_\circ$ to be $(\c\circ(\d\circ\e))\i((\c\circ\d)\circ\e)$.
(See Remark \ref{rem:circinv}.)  Main Theorem \ref{mainthm:isotopy}
then comes from the following result.

\begin{thm}\label{thm:translate}
We have the formulas 
\begin{align}
\c^{\circ n} &= \c^n, \label{eq:transsigma} \\
[\c,\d]_\circ &= \chi(c,d)\alpha(c,k,d)\i, \label{eq:transchi} \\
[\c,\d,\e]_\circ &= \alpha(c,d,e). \label{eq:transalpha}
\end{align}
\end{thm}
In particular, $\c^{\circ 3} = \sigma(\c)$.
\begin{proof}
We first verify \eqref{eq:transsigma} by induction.  Applying
\eqref{eq:kappa}, we have
\begin{equation}
\c^{\circ(n+1)} = \c\circ\c^n = \c\c^n\cdot\alpha(c,k,c^n) = \c^{n+1},
\end{equation}
with the last equality following from the symplectic property of
$\alpha$.  \eqref{eq:transsigma} follows.

Next, from \eqref{eq:kappa} and 
\begin{equation}
\d \circ_\kappa \c = \d\c\cdot\alpha(d,k,c),
\end{equation}
it follows that
\begin{equation}
\begin{split}
[\c,\d]_\circ &= (\d\circ\c)\i(\c\circ\d) \\
	&= (\d\c)\i(\c\d) \cdot \alpha(d,k,c)\i\alpha(c,k,d) \\
	&= [\c,\d] \cdot \alpha(c,k,d)^2 \\
	&= \chi(c,d)\alpha(c,k,d)\i,
\end{split}
\end{equation}
with the last equality following because $Z$ has exponent 3.
\eqref{eq:transchi} follows.

Finally, since
\begin{align}
(\c\circ\d)\circ\e &= (\c\d)\e \cdot \alpha(c,k,d) \alpha(cd,k,e) \\
\c\circ(\d\circ\e) &= \c(\d\e) \cdot \alpha(d,k,e) \alpha(c,k,de),
\end{align}
we have that
\begin{equation}
\begin{split}
[\c,\d,\e]_\circ &= (\c\circ(\d\circ\e))\i((\c\circ\d)\circ\e) \\
	&= [\c,\d,\e]\cdot\alpha(c,k,d)\alpha(cd,k,e)
		\alpha(d,k,e)\i\alpha(c,k,de)\i \\
	&= \alpha(c,d,e)\alpha(c,k,d)\alpha(c,k,e)\alpha(d,k,e) \\
	&\qquad\cdot
		\alpha(d,k,e)\i\alpha(c,k,d)\i\alpha(c,k,e)\i \\
	&= \alpha(c,d,e).
\end{split}
\end{equation}
The theorem follows.
\end{proof}

\begin{proof}[{\bf Proof of Main Theorem \ref{mainthm:isotopy}}]
From \eqref{eq:transsigma}--\eqref{eq:transalpha}, it follows that
$(L,\circ_\kappa)$ is the coded extension of $\adt_{k\i}(C)$, since
condition \ref{cond:centralext} of Definition \ref{defn:codedext}
follows from \eqref{eq:transchi} and \eqref{eq:transalpha}, and the
formulas \eqref{eq:CEpower}--\eqref{eq:CEassociate} of condition
\ref{cond:formulas} of Definition \ref{defn:codedext} follow from the
formulas \eqref{eq:transsigma}--\eqref{eq:transalpha}, respectively.
Therefore, any $\kappa$-isotope of $L$ is a coded extension of an
adjoint translate of $C$, and vice versa.  The theorem follows.
\end{proof}

We illustrate Main Theorems \ref{mainthm:CVSisSFML} and
\ref{mainthm:isotopy} by enumerating the isomorphism and isotopy
classes of the nonassociative coded extensions of the CVS's of
dimension 3 and 4 over $\F_3$.  For simplicity, we only discuss the
exponent 3 ($\sigma=1$) cases.  In the following, we let
$(C,1,\chi,\alpha)$ be the CVS under discussion, we let $L$ be its
coded extension, and we let $Z=\gen{\omega}$.

\begin{exmp}
In the dimension 3 case, either $\chi=1$ (i.e., $L$ is commutative),
or there exists a basis $\set{k,c,d}$ of $C$ such that
$\rad(\chi)=\gen{k}$ and $\chi(c,d)=\omega$.  After possibly inverting
$k$, we may assume that $\alpha(c,k,d)=\omega$, so there are two
possibilities for $L$, up to isomorphism.  However, in the
non-commutative case, $\adt_k(C)$ has $\chi=1$ as its bilinear form,
so the two isomorphism classes are isotopic.
\end{exmp}

\begin{exmp}
In dimension 4, since we assume $\alpha$ is nontrivial, there is some
$c\in C$ such that $c\notin\rad(\alpha)$.  Therefore, since
$\alpha(c,-,-)$ is a nontrivial bilinear symplectic form on $C$ whose
radical contains $c$, we may choose a basis $\set{c,d,e,f}$ for $C$
such that $\alpha(c,d,-)=1$ and $\alpha(c,e,f)=\alpha(d,e,f)\i$.  It
may then be easily verified on the basis $\set{c,d,e,f}$ that
$cd\in\rad(\alpha)$.  Therefore, $\rad(\alpha)$ is 1-dimensional.  It
follows easily that $\alpha$ is unique up to isomorphism, and that we
have four possible isomorphism classes for $C$ (and therefore, for
$L$):
\begin{enumerate}
\item\label{enum:chitriv} $\chi$ trivial;
\item\label{enum:chinon} $\chi$ nondegenerate;
\item\label{enum:radin} $\rad(\chi)$ 2-dimensional, containing
$\rad(\alpha)$; and
\item\label{enum:radnot} $\rad(\chi)$ 2-dimensional, not containing
$\rad(\alpha)$.
\end{enumerate}

So now let $\rad(\alpha)=\gen{c}$.  First, for $d\in C$, $\chi(c,d)$
is invariant under adjoint translation, so isomorphism classes
\ref{enum:chitriv} and \ref{enum:chinon} cannot be adjoint translates.
Conversely, if $\chi$ is nondegenerate, without loss of generality, we
may choose a basis $\set{c,k,d,e}$ for $C$ such that $\gen{c,k}$ and
$\gen{d,e}$ are orthogonal with respect to $\chi$ and
$\chi(d,e)=\alpha(d,k,e)$, in which case $\adt_k(C)$ is in isomorphism
class \ref{enum:radnot}.  On the other hand, if $\rad(\alpha)=\gen{c}$
and $\rad(\chi)=\gen{c,k}$, then by inverting $k$ if necessary, we may
choose a basis $\set{c,k,d,e}$ for $C$ such that
$\chi(d,e)=\alpha(d,k,e)$, in which case $\adt_k(C)$ is in isomorphism
class \ref{enum:chitriv}.  Therefore, we have precisely two isotopy
classes: isomorphism classes \ref{enum:chitriv} and \ref{enum:radin},
and isomorphism classes \ref{enum:chinon} and \ref{enum:radnot}.
\end{exmp}

Compare B\'{e}n\'{e}teau~\cite[IV.3]{lb:cml} and Ray-Chaudhuri and
Roth~\cite{dkrcrr:hallcml}.

\section{A construction of all finite Moufang loops of class 2}
\label{sect:constructclasstwo}

In this section, we give a construction of all finite Moufang
$p$-loops of class 2, which, because of Theorem
\ref{thm:torsionsubloop}, gives a construction of all finite Moufang
loops of class 2.  We will mostly be imitating Sections
\ref{sect:eml}--\ref{sect:isotopy}, using modules instead of
vector spaces, so many details will be omittted.

\begin{notn}
For the rest of this section, we let $p$ be a prime, we let $L$ be a
Moufang $p$-loop with a central subgroup $Z$ such that $C\isom L/Z$ is
an abelian group, and we resume our $\c,\d,\e,\dots$ and $c,d,e,\dots$
convention.
\end{notn}

We would like to say that $L$ gives $C$ the structure of a ``coded
module,'' whatever that means.  Now, $\chi$ and $\alpha$ work exactly
as they do for SFML's (Theorem \ref{thm:sympmulti}).  However, if $r$
is the exponent of $Z$, the only well-defined ``power functions''
$\sigma: C \rightarrow Z$ are the $q$th power functions, where $r$
divides $q$, so if $r>p$, these power functions no longer capture all
of the information we need.  We fix this problem with the following
definition.

\begin{defn}\label{defn:codmod}
A {\it coded module\/} is defined to be a 6-tuple
$(C,Z,\set{c_i},\set{z_i},\chi,\alpha)$, where $C$ and $Z$ are finite
abelian $p$-groups; $\set{c_i}$ is a basis (set of independent
generators) for $C$; $\set{z_i}$ is a set of elements of $Z$; and
$\chi: C\times C\rightarrow Z$ and $\alpha: C\times C\times
C\rightarrow Z$ satisfy \eqref{eq:chisymp}--\eqref{eq:chimultilin} and
\eqref{eq:alphasymp}--\eqref{eq:alphamultilin} for all $c,d,e,f\in C$
and all $n\in\Z$.
\end{defn}

We then see that $L$, along with a choice of basis $\set{c_i}$ for $C$
and a choice of preimages $\c_i$ in $L$, makes $C$ and $Z$ into a
coded module, by defining $z_i = \c_i^{q_i}$, where $q_i$ is the order
of $c_i$, and defining $\chi$ and $\alpha$ as usual.

We may now recover much of Sections \ref{sect:eml} and
\ref{sect:codedextensions} in the case of coded modules.  First, we
can generalize the definition of coded extension (Definition
\ref{defn:codedext}) to coded modules by removing the condition that
$Z$ have order $p$, and replacing \eqref{eq:CEpower} with
\begin{equation}\label{eq:codedmodpower}
\c_i^{q_i} = z_i.
\end{equation}
Imitating Theorem \ref{thm:codedextisSFML}, we see that a loop $L$
(with choice of preimages, etc.) is a coded extension of $C$ if and
only if $L$ is a Moufang $p$-loop of class 2 whose coded module is
$C$.

Having carried over the basic definitions, we next generalize the
results of Section \ref{sect:codedextensions} appropriately.  First,
the uniqueness of coded extensions follows from the proof of Theorem
\ref{thm:codedextunique}, replacing \eqref{eq:Cpowerreln} with
\eqref{eq:codedmodpower}, and \eqref{eq:zpower} with defining (group)
relations for $Z$.  The definition of the semidirect central product
(Definition \ref{defn:sdcp}) also still works, as does Theorem
\ref{thm:sdcpiscodedext}, if we replace the verification of
\eqref{eq:CEpower} by the observation that concatenation of the bases
of $D$ and $E$ produces a basis for $C$.  Finally, coded extensions of
1-dimensional coded modules are easy to construct using central
products of groups, so Theorem \ref{thm:codedextexist} follows as
before.

In fact, in general, the only result not carried over from
Section~\ref{sect:codedextensions} is the ``isomorphic SFML's implies
isomorphic CVS's'' statement of Main Theorem \ref{mainthm:CVSisSFML},
since the structure of a coded module is highly basis-dependent.
However, if $Z$ is elementary abelian, then by defining $\sigma_q(c):
C_q\rightarrow Z$ by $\sigma_q(c)=\gamma^q$, for all $q=p^n$, where
$C_q$ is the subgroup of $C$ of all elements of $C$ whose order
divides $q$, we may actually recover all of Main Theorem
\ref{mainthm:CVSisSFML}.  In fact, since \eqref{eq:sigmalin} holds in
the elementary abelian case if we replace $p$ with $q$, everything
works as before.  (In particular, when $q$ is a power of 2, $\sigma_q$
is linear for $q>2$.)

We can also generalize Main Theorem \ref{mainthm:isotopy} to the case
where $Z$ is an elementary abelian 3-group and $C$ is a coded module
with values in $Z$.  We state this result in the following theorem,
whose proof may be copied directly from Section \ref{sect:isotopy}.
(Note that the proof of Theorem \ref{thm:translate} only uses the
exponent of $Z$, and not its order.)

\begin{thm}\label{thm:expthreeisotopy}
Let $L$ be a coded extension of a coded module $C$ of 3-power order by
an elementary abelian 3-group.  Then up to isomorphism, the
loop-isotopes of $L$ are precisely the coded extensions of the adjoint
translates of $C$.  In particular, $\sigma_q$ ($q=3^n$) and $\alpha$
are ``isotopy invariants'' of $L$.  \qed
\end{thm}

For example, Theorem \ref{thm:expthreeisotopy} applies when $L$ has
exponent 3 and class 2.

Finally, we address the question analogous to the one posed in Section
\ref{sect:cvsiscode}: that of which $\chi$, $\alpha$, etc., can be
chosen for a coded module.  (The 1-dimensional case shows that any
values may be chosen for the $z_i$.)  Again, for $p>2$, this is easy
(keeping in mind that $\alpha^6=1$), so it will be enough to prove
Theorem \ref{thm:consistenttwo}, in which we use the following
conventions.

\begin{notn}
We revert to additive notation for abelian groups.  All summations,
indices, etc., are from 1 to $n$ in the manner indicated.  (For
instance, $i<j$ means for all $1\le i<j\le n$.)  We also let
$\alpha_{ijk}=\alpha(x_i,x_j,x_k)$.
\end{notn}

\begin{thm}\label{thm:consistenttwo}
Let $C$ and $Z$ be finite abelian 2-groups, and let $\set{x_i}$ be a
basis for $C$.  Choose some symplectic multilinear function $\alpha:
C\times C\times C \rightarrow Z$ such that $2\alpha = 0$ identically,
and for all $i,j$, choose $\chi_{ij}\in Z$ such that:
\begin{enumerate}
\item\label{chicond:symp} $\chi_{ii}=0$,
\item\label{chicond:skew} $\chi_{ij}=-\chi_{ji}$, and
\item\label{chicond:order} the (additive) order of $\chi_{ij}$ divides
the order of $x_i$.
\end{enumerate}
Let $c=\sum_i c_i x_i$ and $d=\sum_j d_j x_j$.  The function $\chi:
C\times C\rightarrow Z$ defined by
\begin{equation}\label{eq:chilongformula}
\chi(c,d) =
	\sum_{i\neq j} c_i d_j \chi_{ij}
	+ \sum_{i<j} \sum_k c_i c_j d_k \alpha_{ijk}
	+ \sum_{i} \sum_{j<k} c_i d_j d_k \alpha_{ijk},
\end{equation}
along with any choice of elements $z_i\in Z$, gives $C$ the structure
of a coded module.  Furthermore, any even coded module may be
constructed in this way.
\end{thm}
Note that the first term of the right-hand side of
\eqref{eq:chilongformula} makes sense only when condition
\ref{chicond:order} is satisfied, since the $c_i$ in $c=\sum_i c_i
x_i$ is only well-defined as an integer mod the order of $x_i$.  Also,
the last two terms of \eqref{eq:chilongformula} are invariant under a
reordering of the basis only because $\alpha_{ijk}$ is symmetric in
all three indices.
\begin{proof}
Given a coded module, if $\chi_{ij}=\chi(x_i,x_j)$, then conditions
\ref{chicond:symp}, \ref{chicond:skew}, and \ref{chicond:order} follow
from \eqref{eq:chisymp}, \eqref{eq:chiskew}, and Theorem
\ref{thm:gcd}, respectively.  Furthermore, \eqref{eq:chilongformula} follows
from repeated application of \eqref{eq:chimultilin} and
\eqref{eq:alphamultilin}, and $2\alpha=0$ follows from Main Theorem
\ref{mainthm:assocsix}.  Our last assertion follows easily.

As for the first assertion, it suffices to check
\eqref{eq:chisymp}--\eqref{eq:chimultilin}; in fact, because of
\eqref{eq:chipowerlintrick}, it suffices to check \eqref{eq:chisymp},
\eqref{eq:chiskew}, and \eqref{eq:chimultilin}.  In the following, let
$c=\sum_i c_i x_i$, $d=\sum_i d_i x_i$, and $e=\sum_i e_i x_i$.

To check \eqref{eq:chisymp}, we compute
\begin{equation}
\begin{split}
\chi(c,c) &=
	\sum_{i\neq j} c_i c_j \chi_{ij}
	+ \sum_{i<j} \sum_k c_i c_j c_k \alpha_{ijk}
	+ \sum_{i} \sum_{j<k} c_i c_j c_k \alpha_{ijk} \\
&=
	\sum_{i<j} \(c_i c_j \chi_{ij} + c_j c_i \chi_{ji}\)
	+ \sum_{i<j} \sum_k c_i c_j c_k \alpha_{ijk}
	+ \sum_{j<k} \sum_{i} c_j c_k c_i \alpha_{jki} \\
&=	0,
\end{split}
\end{equation}
since $\chi_{ij}+\chi_{ji}=0$, $\alpha_{jki}=\alpha_{ijk}$, and
$2\alpha=0$.

To check \eqref{eq:chiskew}, we compute
\begin{equation}
\begin{split}
\chi(c,d) &=
	\sum_{i\neq j} c_i d_j \chi_{ij}
	+ \sum_{i<j} \sum_k c_i c_j d_k \alpha_{ijk}
	+ \sum_{i} \sum_{j<k} c_i d_j d_k \alpha_{ijk} \\
&=
	- \sum_{j\neq i} d_j c_i \chi_{ji}
	+ \sum_{j<k} \sum_i d_j d_k c_i \alpha_{jki}
	+ \sum_k \sum_{i<j} d_k c_i c_j \alpha_{kij} \\
&=	- \chi(d,c),
\end{split}
\end{equation}
since $\chi_{ij}=-\chi_{ji}$,
$\alpha_{jki}=\alpha_{kij}=\alpha_{ijk}$, and $\alpha=-\alpha$.

Finally, to check \eqref{eq:chimultilin}, we compute
\begin{equation}
\begin{split}
\chi(c+d,e) &=
	\sum_{i\neq j} (c_i + d_i) e_j \chi_{ij}
	+ \sum_{i<j} \sum_k (c_i + d_i) (c_j + d_j) e_k \alpha_{ijk}
	\\
&\qquad
	+ \sum_{i} \sum_{j<k} (c_i + d_i) e_j e_k \alpha_{ijk}
	\\
&=
	\sum_{i\neq j} c_i e_j \chi_{ij}
	+ \sum_{i\neq j} d_i e_j \chi_{ij} \\
&\qquad
	+ \sum_{i<j} \sum_k c_i c_j e_k \alpha_{ijk}
	+ \sum_{i<j} \sum_k d_i d_j e_k \alpha_{ijk} \\
&\qquad
	+ \sum_{i} \sum_{j<k} c_i e_j e_k \alpha_{ijk}
	+ \sum_{i} \sum_{j<k} d_i e_j e_k \alpha_{ijk} \\
&\qquad
	+ \sum_{i<j} \sum_k c_i d_j e_k \alpha_{ijk}
	+ \sum_{i<j} \sum_k d_i c_j e_k \alpha_{ijk}
	\\
&=	\chi(c,e) + \chi(d,e) + \alpha(c,d,e),
\end{split}
\end{equation}
since
\begin{equation}
\begin{split}
&\sum_{i<j} \sum_k c_i d_j e_k \alpha_{ijk}
	+ \sum_{i<j} \sum_k d_i c_j e_k \alpha_{ijk} \\
{}={}
&\sum_{i<j} \sum_k c_i d_j e_k \alpha_{ijk}
	+ \sum_{i>j} \sum_k c_i d_j e_k \alpha_{ijk} \\
{}={} &\alpha(c,d,e),
\end{split}
\end{equation}
using $\alpha_{ijk}=\alpha_{jik}$ and $\alpha_{kjk}=\alpha_{ikk}=0$.
The theorem follows.
\end{proof}

Finally, when $Z$ is an elementary abelian $p$-group, we have seen
previously that it is possible to define $\sigma_q: C_q\rightarrow Z$
for all $q=p^n$, and so the problem arises of which values can be
chosen for $\sigma_q$.  Now, since $\sigma_q$ is trivial on $C_{q/p}$,
for $q>2$, we may choose $\sigma_q$ freely on a basis for
$C_q/C_{q/p}$, since $\sigma_q$ is linear for $q>2$.  Therefore, the
only case remaining is $\sigma=\sigma_2$.  We leave it to the
interested reader to verify (by imitating the proof of Theorem
\ref{thm:consistenttwo}) that if we choose arbitrary elements
$\sigma_i\in Z$, and define
\begin{equation}\label{eq:sigmatwo}
\sigma(c) = \sum_i c_i \sigma_i + \sum_{i<j} c_i c_j \chi_{ij}
		+ \sum_{i<j<k} c_i c_j c_k \alpha_{ijk},
\end{equation}
then $\sigma$ has the desired properties.

\begin{rem}
Note that since our results, in most cases, say nothing about the
isomorphism problem for Moufang loops of class 2, the result in group
theory which is probably closest to our construction is the fact that
every nilpotent group has a consistent polycyclic presentation.  For
more about polycyclic presentations, see
Sims~\cite[9.4]{ccs:compfpgp}.
\end{rem}

%%%%% end matter

\section{Acknowledgements}

The author was partly supported by a University of Michigan Rackham
Summer Faculty Fellowship, and by the generosity of MSRI and its
staff.  (Research at MSRI is supported in part by NSF grant
DMS-9022140.)  The author would also like to thank J. H. Conway and
L. Schneps for their helpful remarks.  Extraspecial thanks goes to
R. L. Griess for making many helpful comments on an early draft of
this paper, and for considerable help in formulating Definition
\ref{defn:SFL} and proving Theorem \ref{thm:SFMLisCSFML}.

%%\bibliographystyle{amsplain}
%%\bibliography{Bib/strings,Bib/groups,Bib/knots,Bib/geom,Bib/pubs,%
%%Bib/automorphic,Bib/algebra,Bib/loops}
\providecommand{\bysame}{\leavevmode\hbox to3em{\hrulefill}\thinspace}

\end{document}